\documentclass[11pt, a4paper]{article}

\usepackage[utf8]{inputenc}
\usepackage[T1]{fontenc}
\usepackage{amsmath, amssymb, amsthm}
\usepackage{geometry}
\usepackage{tikz}
\usetikzlibrary{arrows.meta, positioning, shapes.geometric, calc}
\usepackage{enumitem}
\usepackage{hyperref}
\usepackage{cleveref}

\newtheorem{theorem}{Theorem}[section]
\newtheorem{lemma}[theorem]{Lemma}
\newtheorem{proposition}[theorem]{Proposition}
\newtheorem{corollary}[theorem]{Corollary}
\newtheorem{conjecture}[theorem]{Conjecture}
\theoremstyle{definition}
\newtheorem{definition}[theorem]{Definition}
\newtheorem{example}[theorem]{Example}
\newtheorem{remark}[theorem]{Remark}
\newtheorem{observation}[theorem]{Observation}

\crefname{theorem}{theorem}{theorems}
\crefname{lemma}{lemma}{lemmas}
\crefname{proposition}{proposition}{propositions}
\crefname{corollary}{corollary}{corollaries}
\crefname{conjecture}{conjecture}{conjectures}
\crefname{definition}{definition}{definitions}
\crefname{example}{example}{examples}
\crefname{remark}{remark}{remarks}
\crefname{observation}{observation}{observations}
\crefname{question}{question}{questions}

\Crefname{theorem}{Theorem}{Theorems}
\Crefname{lemma}{Lemma}{Lemmas}
\Crefname{proposition}{Proposition}{Propositions}
\Crefname{corollary}{Corollary}{Corollaries}
\Crefname{conjecture}{Conjecture}{Conjectures}
\Crefname{definition}{Definition}{Definitions}
\Crefname{example}{Example}{Examples}
\Crefname{remark}{Remark}{Remarks}
\Crefname{observation}{Observation}{Observations}
\Crefname{question}{Question}{Questions}

\hypersetup{
  pdfencoding=auto,
  psdextra,
  pdftitle={Hammock coordinates for silting subcategories of derived-discrete algebras},
  pdfauthor={Shantanu Sardar},
}

\newcommand{\Db}{\mathcal{D}^b}

\title{Hammock coordinates for silting subcategories of derived-discrete algebras}
\author{Shantanu Sardar\thanks{Shaanxi Normal University, China;
  \texttt{shantanusardar17@gmail.com}}}
\date{}

\begin{document}

\maketitle

\begin{abstract}
We introduce hammock coordinates to study the combinatorial structure of partial silting configurations in the bounded derived category of a derived-discrete algebra $\Lambda(r,n,m)$. By extracting exact integer inequalities from the Bobi\'nski--Gei\ss{}--Skowro\'nski Hom-hammock formula, we establish a closed-form criterion for pairwise partial silting both within any $\mathcal Z$-component and across distinct components. Using this coordinate framework, we establish three primary results. First, we prove the universal upper bound $|S| \le n+m$ on the cardinality of any partial silting configuration $S$ of $\mathcal Z$-objects
---both globally via Bongartz completion in the silting-discrete setting, and combinatorially within a single component via a $\Sigma^r$-invariant staircase argument. Second, we prove that this bound is universally sharp by constructing an explicit L-shaped partial silting configuration of size $n+m$ for every $\Lambda(r,n,m)$, alongside specialized diagonal configurations. Third, we show that the configuration lattice carries a faithful action of the cyclic group $\mathbb Z/g$ where $g = \gcd(r+m, n-r)$. We conclude by showing a structural obstruction to local determination: pairwise partial silting of consecutive objects in the lattice ordering does not imply global partial silting, due to the strict additivity of a coordinate defect along chains.


\end{abstract}

\section{Introduction}
\label{sec:intro}

In the representation-theoretic study of finite dimensional associative algebras over an algebraically closed field $K$, Brenner introduced certain partially ordered sets known as hammocks to encode factorisations of maps between indecomposable finitely generated modules \cite{Brenner1986}. The classification of finite dimensional $K$-algebras whose bounded derived category contains only finitely many isomorphism classes of indecomposable complexes of each fixed homology dimension vector was initiated by Vossieck \cite{Vossieck2001} and completed by Bobi\'nski, Gei\ss{}, and Skowro\'nski \cite{BGS2004}; these are the derived-discrete algebras. The Auslander--Reiten (A-R) quiver---encoding the isomorphism classes of indecomposable objects and irreducible morphisms between them---was explicitly described in \cite{BGS2004}. It consists of certain $\mathcal{X}$- and $\mathcal{Y}$-components, together with $r$ components of type $\mathbb{Z}A_\infty^\infty$, denoted $\mathcal{Z}^0, \dots, \mathcal{Z}^{r-1}$. Crucially, \cite{BGS2004} provides a Hom-hammock formula characterising the non-vanishing of morphisms between indecomposable $\mathcal{Z}$-objects.

Silting theory, introduced by Aihara and Iyama \cite{AI2012}, provides a framework generalising tilting theory. Silting subcategories, which are presilting subcategories that generate the bounded derived category as a thick subcategory, control bounded $t$-structures and co-$t$-structures. In the context of derived-discrete algebras, Broomhead, Pauksztello, and Ploog classified silting objects up to shift in \cite{BPP2017}, and showed in \cite{BPP2016} that the associated silting-pairs CW complex is contractible.

However, the geometric properties of the silting poset in \cite{BPP2016, BPP2017} rest on abstract arguments involving length categories. The translation of the pairwise partial silting condition into explicit, computable combinatorial data had not previously been carried out in the literature. The absence of a closed-form combinatorial criterion limits the ability to explicitly construct maximal configurations or to test local-to-global properties within the silting poset. The primary objective of the present paper is to address this gap by translating the homological vanishing conditions into exact integer inequalities on a two-dimensional lattice.

\subsection{Main results}
\label{subsec:intro-main}

The first major achievement of the paper is the development of a closed-form criterion for pairwise partial silting. We extract exact integer inequalities from the Hom-hammock formula to completely characterise when $\operatorname{Hom}^{>0}(Z_1, Z_2) = 0$ for any two objects in the $\mathcal Z$-components of $\Db(\Lambda(r,n,m))$, reducing homological computations to lattice combinatorics.

\begin{theorem}[see \Cref{thm:pairwise} and \Cref{thm:cross-pairwise}]\label{thm:intro-criterion}
Let $\Lambda = \Lambda(r,n,m)$ with $n > r \ge 1$ and $m \ge 0$, and set $a = r+m$ and $b = n-r$. Let $Z_1 = Z_{i_1,j_1}^{k_1} \in \mathcal Z^{k_1}$ and $Z_2 = Z_{i_2,j_2}^{k_2} \in \mathcal Z^{k_2}$ be distinct indecomposable objects, and set $U = i_1 - i_2$ and $V = j_2 - j_1$. If $k_1 = k_2$, orient $Z_1, Z_2$ so that $V \ge 0$. Then $\{Z_1, Z_2\}$ is a partial silting pair if and only if one of the following holds:
\begin{enumerate}[label=\textup{(\arabic*)}]
  \item[\textup{(A)}] $V = 0$ and $1 \le |U| \le a$;
  \item[\textup{(B1)}] $1 \le V \le b - 1$ and $0 \le U \le a$;
  \item[\textup{(B2)}] $V \ge b$ and $\lfloor (U-1)/a \rfloor = \lfloor V/b \rfloor$.
\end{enumerate}
If $0 \le k_1 < k_2 \le r-1$, then $\{Z_1, Z_2\}$ is a partial silting pair if and only if exactly one of the following holds: \textup{(i)} $V = 0$ and $-a \le U \le -1$; \textup{(ii)} $V \ge 1$, $U \ge 0$ satisfy \textup{(B1)} or \textup{(B2)}; \textup{(iii)} $-V \ge 1$, $-U \ge 1$ satisfy \textup{(B1)} or \textup{(B2)} evaluated on $(-U, -V)$.
\end{theorem}

The second and perhaps the most important achievement is the establishment of the universal upper bound $n+m$ on the size of any partial silting configuration of $\mathcal Z$-objects, and the proof that this bound is universally sharp.

\begin{theorem}[see \Cref{thm:sharp-combinatorial} and \Cref{thm:bound}]\label{thm:intro-bound}
Let $S$ be a partial silting configuration of pairwise non-isomorphic $\mathcal Z$-objects in $\Db(\Lambda(r,n,m))$. Then $|S| \le n+m$. This bound holds combinatorially within a single component via a monotonic $\Sigma^r$-invariant $(I,r)$-staircase, and globally across all components via Bongartz completion.
\end{theorem}

\begin{theorem}[see \Cref{thm:sharp-universal}]\label{thm:intro-sharp}
The upper bound $n+m$ is universally sharp. The explicit L-shaped configuration
\[
  S \;:=\; \bigl\{\, Z_{-p,\,0}^0 \;\big|\; 0 \le p \le r+m-1 \,\bigr\}
  \;\sqcup\;
  \bigl\{\, Z_{-(r+m),\,q}^0 \;\big|\; 0 \le q \le n-r-1 \,\bigr\}
\]
is a partial silting configuration of size $n+m$ for every $\Lambda(r,n,m)$.
\end{theorem}

A less significant contribution of the paper, explored here for various reasons, is the structural analysis of the lattice symmetries and obstructions, summarized in the final theorem.

\begin{theorem}[see \Cref{thm:Zg-action} and \Cref{thm:notlocal}]\label{thm:intro-structure}
The configuration lattice carries a faithful action of the cyclic group $\mathbb Z/g$, where $g = \gcd(r+m, n-r)$. Furthermore, local determination fails: pairwise partial silting of consecutive objects in the $j$-ordering does not imply global partial silting, due to the strict additivity of a coordinate defect $\partial(Z_1, Z_2) = bU - aV$ along chains.
\end{theorem}

\subsection{Flowchart of the coordinate reduction}
\label{subsec:intro-flowchart}

The following flowchart succinctly depicts the steps of the combinatorial reduction used to compute the partial silting criteria and construct sharp configurations.

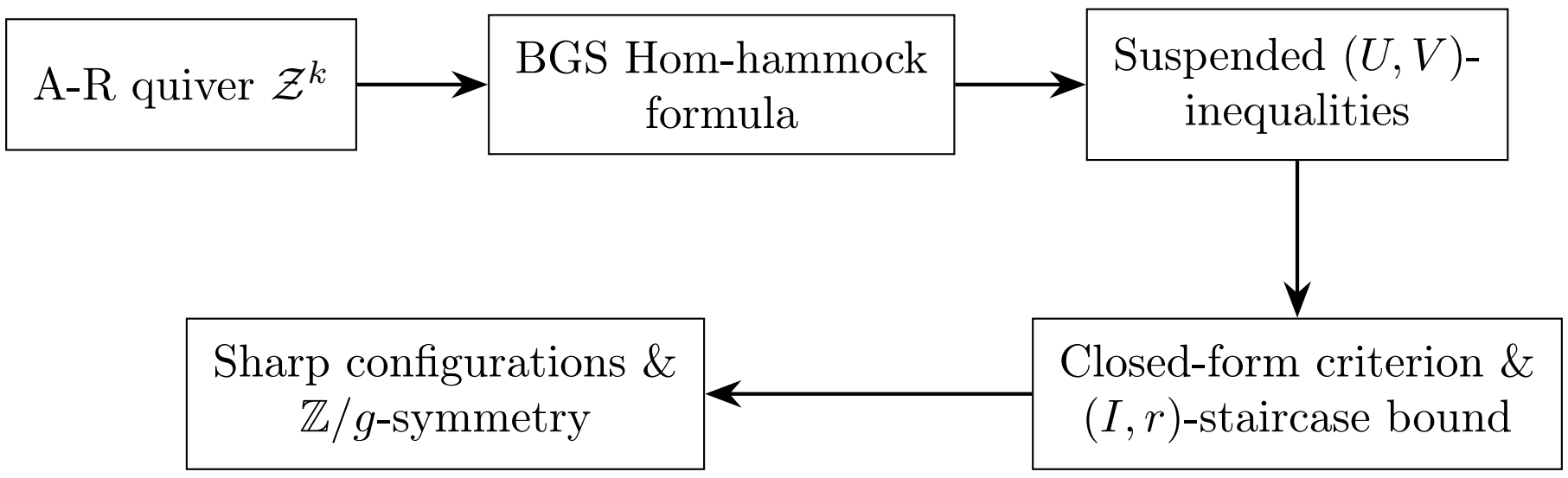
\begin{center}
\begin{tikzpicture}[
  node distance=1.2cm and 1.0cm,
  box/.style={draw, rectangle, align=center, inner sep=2mm, minimum height=1cm, font=\small},
  arrow/.style={-{Stealth[scale=1.2]}, thick}
]

\node[box] (ar) {A-R quiver $\mathcal{Z}^k$};
\node[box, right=of ar] (bgs) {BGS Hom-hammock\\formula};
\node[box, right=of bgs] (ineq) {Suspended $(U,V)$-\\inequalities};
\node[box, below=of ineq] (closed) {Closed-form criterion \&\\$(I,r)$-staircase bound};
\node[box, left=of closed, xshift=-1.5cm] (sharp) {Sharp configurations \&\\$\mathbb{Z}/g$-symmetry};

\draw[arrow] (ar) -- (bgs);
\draw[arrow] (bgs) -- (ineq);
\draw[arrow] (ineq) -- (closed);
\draw[arrow] (closed) -- (sharp);

\end{tikzpicture}
\end{center}

\subsection{Contents of sections}
\label{subsec:intro-contents}

Due to the highly complex nature of the combinatorics of string algebras and derived categories, the paper is very heavy on notation and terminology. To facilitate reading we have included several remarks as well as (running) examples to illustrate new terminology. The paper is organised as follows.

Starting with the definition of the derived-discrete algebra $\Lambda(r,n,m)$, we recall fundamental concepts like the Auslander--Reiten quiver, the Hom-hammock formula, and silting theory in \S 2---the treatment mostly follows the established literature. The first major achievement, namely the translation of the suspended hammock conditions into $(U,V)$-coordinates and the closed-form intra- and cross-component pairwise partial silting criteria, is detailed in \S 3.

To quantify the conditions obtained in \S 3, we introduce the $L$-invariant and the coordinate defect $\partial$ in \S 4, and study their arithmetic properties along chains. The second step in the computation of maximal configurations is the universal upper bound $k \le n+m$, which is proved in \S 5. We provide both the Bongartz completion argument and the unconditional combinatorial $(I,r)$-staircase proof.

The next key result shows that the $n+m$ bound is universally sharp. In \S 6, we construct the universal L-shaped configuration for arbitrary $\Lambda(r,n,m)$ and the diagonal chain for $\Lambda(1,n,0)$. Finally, to complete the structural analysis of the lattice, we establish the faithful $\mathbb{Z}/g$-action on partial silting configurations in \S 7, and utilise the defect arithmetic in \S 8 to illustrate the failure of local determination.

We conclude the paper by comparing our hammock data with the Avella-Alaminos--Gei\ss{} derived invariant and the Hom-hammocks of Broomhead, Pauksztello, and Ploog in \S 9. Future directions and related open problems are discussed in \S 10. Appendix \ref{app:tables} contains explicit case-by-case verifications of the pairwise criteria, and Appendix \ref{app:algorithm} presents a clique enumeration algorithm in hammock coordinates.

\section{Preliminaries on derived-discrete algebras and silting}
\label{sec:preliminaries}

We begin by introducing the canonical family of derived-discrete algebras $\Lambda(r,n,m)$. After describing its Auslander--Reiten (A-R) quiver, we recall the Hom-hammock formula characterising morphisms between $\mathcal{Z}$-objects. In the latter half, we review essential results regarding presilting and silting subcategories and their connections to algebraic $t$-structures via the K\"onig--Yang correspondence.

\subsection{The canonical family \texorpdfstring{$\Lambda(r,n,m)$}{Lambda(r,n,m)} and its Auslander--Reiten quiver}\label{subsec:prelim-AR}

Let $K$ be an algebraically closed field. The derived-discrete algebras are finite dimensional $K$-algebras whose bounded derived category contains only finitely many isomorphism classes of indecomposable complexes of each fixed homology dimension vector. A representative class is given by the family $\Lambda(r,n,m)$.

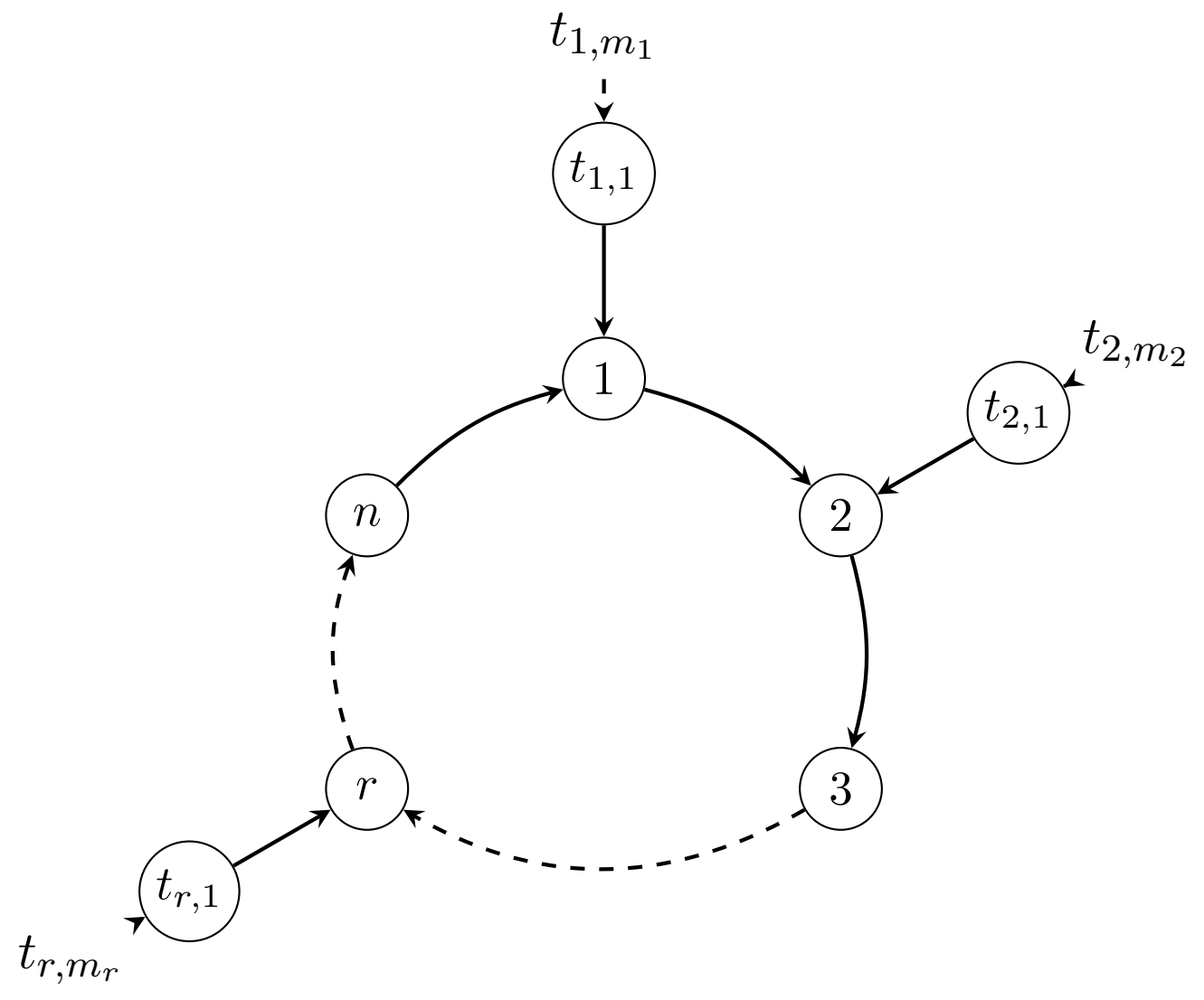
\begin{figure}[htbp]
\centering
\begin{tikzpicture}[
  >=stealth,
  node distance=1.5cm,
  vertex/.style={circle, draw, inner sep=1.5pt, minimum size=6mm, font=\small}
]
  \node[vertex] (v1) at (90:2) {$1$};
  \node[vertex] (v2) at (30:2) {$2$};
  \node[vertex] (v3) at (330:2) {$3$};
  \node[vertex] (vr) at (210:2) {$r$};
  \node[vertex] (vn) at (150:2) {$n$};

  \draw[->, thick] (v1) to[bend left=15] (v2);
  \draw[->, thick] (v2) to[bend left=15] (v3);
  \draw[->, thick, dashed] (v3) to[bend left=30] (vr);
  \draw[->, thick, dashed] (vr) to[bend left=20] (vn);
  \draw[->, thick] (vn) to[bend left=15] (v1);

  \node[vertex] (t11) at (90:3.5) {$t_{1,1}$};
  \node (t1m) at (90:4.5) {$t_{1, m_1}$};
  \draw[->, thick] (t11) -- (v1);
  \draw[->, thick, dashed] (t1m) -- (t11);

  \node[vertex] (t21) at (30:3.5) {$t_{2,1}$};
  \node (t2m) at (30:4.5) {$t_{2, m_2}$};
  \draw[->, thick] (t21) -- (v2);
  \draw[->, thick, dashed] (t2m) -- (t21);

  \node[vertex] (tr1) at (210:3.5) {$t_{r,1}$};
  \node (trm) at (210:4.5) {$t_{r, m_r}$};
  \draw[->, thick] (tr1) -- (vr);
  \draw[->, thick, dashed] (trm) -- (tr1);
\end{tikzpicture}
\caption{The typical quiver $Q(r,n,m)$ for the derived-discrete algebra $\Lambda(r,n,m)$, featuring a cycle of $n$ vertices and $r$ oriented tails whose lengths $m_1, \dots, m_r$ sum to $m$.}
\label{fig:quiver-rnm}
\end{figure}

For integers $r \ge 1$, $n > r$, and $m \ge 0$, the canonical algebra $\Lambda(r,n,m)$ of \cite{BGS2004} is presented by a bound quiver $(Q(r,n,m), \rho)$ consisting of an oriented cycle of $n$ vertices with $r$ consecutive zero-relations of length two along the cycle, and a single linear tail of length $m$ attached at the cycle. (More generally, distributing the $m$ arrows into $r$ linear tails of lengths $m_1, \dots, m_r$ with $\sum_{k=1}^r m_k = m$ attached to the zero-relation vertices yields a derived-equivalent algebra; see \Cref{fig:quiver-rnm}). In all cases, the total number of vertices of $Q(r,n,m)$ is $n+m$, which exactly matches the rank of $K_0(\Db(\Lambda(r,n,m)))$.

\begin{example}[Running examples]\label{ex:running}
We consider four specific algebras from this family which will serve as running examples to illustrate the hammock combinatorics.
\begin{enumerate}
    \item $\Lambda(1,2,0)$: the Kronecker-like gentle algebra with a $2$-cycle, one zero-relation, and no tails.
    \item $\Lambda(1,3,0)$: the algebra with a $3$-cycle, one zero-relation, and no tails.
    \item $\Lambda(1,3,1)$: the algebra with a $3$-cycle and a single tail of length $1$ attached to the first vertex.
    \item $\Lambda(2,3,0)$: the algebra with a $3$-cycle, two zero-relations, and no tails. \hfill $\triangle$
\end{enumerate}
\end{example}

The bounded derived category $\Db(\Lambda(r,n,m))$ is a triangulated category with suspension functor $\Sigma=[1]$. Because $\Lambda(r,n,m)$ is gentle, it is Gorenstein, and standard Auslander--Reiten theory applies.

\begin{theorem}[{\cite[Theorem~B]{BGS2004}}]\label{thm:AR-structure}
The Auslander--Reiten quiver of $\Db(\Lambda(r,n,m))$ has exactly $3r$ components: $2r$ components of type $\mathbb{Z}A_\infty$ (the $\mathcal{X}$- and $\mathcal{Y}$-components), and $r$ components of type $\mathbb{Z}A_\infty^\infty$ (the $\mathcal{Z}$-components $\mathcal{Z}^0,\dots,\mathcal{Z}^{r-1}$). The indecomposable objects of $\mathcal{Z}^k$ are parameterised as $Z_{i,j}^k$ for $i,j\in\mathbb{Z}$, with irreducible morphisms $Z_{i,j}^k\to Z_{i+1,j}^k$ and $Z_{i,j}^k\to Z_{i,j+1}^k$.
\end{theorem}

Every basic silting object of $\Db(\Lambda(r,n,m))$ has exactly $n+m$ summands \cite{BPP2017}, and at least one of them lies in a $\mathcal{Z}$-component. Indeed, basic silting objects form a $K_0$-basis (\Cref{lem:silting-basic}) of $K_0(\Db(\Lambda))$, and the classes of the $\mathcal{X}$- and $\mathcal{Y}$-components span a proper subgroup; hence any basis must contain at least one $\mathcal{Z}$-summand. The suspension acts on the $\mathcal{Z}$-components by
\begin{equation}\label{eq:sigma-action}
\Sigma\,Z_{i,j}^k = Z_{i,j}^{k+1}\quad (0\le k<r-1),\qquad
\Sigma\,Z_{i,j}^{r-1} = Z_{i+r+m,\,j+r-n}^0.
\end{equation}

\subsection{The Hom-hammock formula}
\label{subsec:prelim-hammock}

The following formula governs the non-vanishing of morphisms between $\mathcal{Z}$-objects. It is the cornerstone for the coordinate criteria developed in this paper.

\begin{theorem}[Hom-hammock formula, {\cite[Section~2]{BGS2004}}]\label{thm:hammock}
Let $A=Z_{i,j}^k\in\operatorname{ind}(\mathcal Z^k)$.
\begin{enumerate}[label=\textup{(\roman*)}]
\item For $B\in\operatorname{ind}(\mathcal Z^k)$, one has $\operatorname{Hom}(A,B)\neq 0$ if and only if $B\cong Z_{a,b}^k$ for some $a\ge i$ and $b\ge j$.
\item For $B\in\operatorname{ind}(\mathcal Z^{k+1})$ (where $\mathcal Z^{k+1}=\mathcal Z^0$ if $k=r-1$), one has $\operatorname{Hom}(A,B)\neq 0$ if and only if $B\cong \Sigma Z_{a,b}^k$ for some $a\le i-1$ and $b\le j-1$.
\item For $B\in\operatorname{ind}(\mathcal Z^{k'})$ with $k'\notin\{k,(k+1)\bmod r\}$, one has $\operatorname{Hom}(A,B)=0$.
\end{enumerate}
When $k = r-1$, the coordinate systems of $\mathcal{Z}^{r-1}$ and $\mathcal{Z}^0$ differ by the shift $(r+m, r-n)$; the condition on $(a,b)$ in \textup{(ii)} is therefore a condition on the lifted coordinates, not on the standard coordinates of $\mathcal{Z}^0$.
\end{theorem}

\begin{remark}\label{rem:hammock-arbitrary-tails}
In \cite{BGS2004}, the canonical algebra $\Lambda(r,n,m)$ is defined with a single tail of length $m$ (corresponding to $m_1 = m$ and $m_2 = \dots = m_r = 0$). By the derived invariance of the Avella-Alaminos--Gei\ss{} invariant \cite{AAG2008}, distributing the $m$ tail arrows across the $r$ zero-relation vertices with lengths $m_1 + \cdots + m_r = m$ yields a derived-equivalent gentle algebra whose $\mathcal{Z}^k$-components, $\Sigma^r$-shift \eqref{eq:sigma-action}, and Hom-hammock formula (\Cref{thm:hammock}) depend only on the sum $m = \sum_{k=1}^r m_k$.
\end{remark}

\subsection{Silting and partial silting subcategories}
\label{subsec:prelim-silting}

Let $\mathcal{T}$ be a $K$-linear, Hom-finite, Krull--Schmidt triangulated category, such as $\Db(\Lambda(r,n,m))$.

\begin{definition}\label{def:silting}
An object $M \in \mathcal{T}$ is called \emph{presilting} if $\operatorname{Hom}_{\mathcal{T}}^{>0}(M, M) = 0$, and \emph{partial silting} if it is a direct summand of a silting object (which in a silting-discrete category is equivalent to being presilting by \Cref{thm:AM}). An object $M$ is called \emph{silting} if it is presilting and $\operatorname{thick}(M) = \mathcal{T}$.
\end{definition}

When dealing with basic objects, we use the terms subcategory and object interchangeably by identifying a basic object with the subcategory $\operatorname{add}(M)$ formed by its direct summands. If $S = \{X_1, \ldots, X_k\}$ is a set of pairwise non-isomorphic indecomposable objects, we say that $S$ is a \emph{partial silting configuration} if $\operatorname{Hom}^{>0}(X_p, X_q) = 0$ for all $1 \le p, q \le k$.

\begin{lemma}[$K_0$-basis for basic silting objects, {\cite[Theorem~2.27]{AI2012}}]\label{lem:silting-basic}
Let $\mathcal T$ be a Hom-finite Krull--Schmidt triangulated category with a basic silting object $T$. Then the isomorphism classes of the indecomposable summands of $T$ form a $\mathbb Z$-basis of $K_0(\mathcal T)$. In particular, every basic silting object of $\mathcal T$ has exactly $\operatorname{rk}K_0(\mathcal T)$ non-isomorphic indecomposable summands.
\end{lemma}

Say that a triangulated category $\mathcal{T}$ is silting-discrete if every interval in the silting poset contains only finitely many basic silting objects; equivalently, the Hasse quiver of the silting poset is locally finite.

\begin{theorem}[{\cite[Appendix~A]{AH2024}}]\label{thm:AH}
Let $\Lambda$ be a derived-discrete algebra. Then the bounded derived category $\Db(\Lambda)$ is silting-discrete.
\end{theorem}

Aihara and Mizuno proved that Bongartz completion---a classical result from $\tau$-tilting theory---generalises to triangulated categories precisely when the category is silting-discrete.

\begin{theorem}[{\cite[Theorem~2.15]{AM2017}}]\label{thm:AM}
Let $\mathcal T$ be a silting-discrete triangulated category and let $M$ be a presilting object of $\mathcal T$. Then $M$ can be completed to a basic silting object $T = M\oplus M'$ of $\mathcal T$.
\end{theorem}

\subsection{The K\"onig--Yang correspondences}
\label{subsec:prelim-KY}

The geometric and combinatorial importance of silting subcategories stems from their exact correspondence with bounded $t$-structures. K\"onig and Yang formalised this bijection.

\begin{theorem}[{\cite[Theorem~0.1]{KY2014}; see also \cite[Section~6]{BPP2017}}]\label{thm:KY}
Let $\Lambda$ be a finite dimensional $K$-algebra. There are canonical bijections between the following sets:
\begin{enumerate}[label=\textup{(\roman*)}]
\item equivalence classes of basic silting objects in $K^b(\operatorname{proj}\Lambda)$;
\item bounded $t$-structures on $\Db(\Lambda)$ whose heart is a length category;
\item bounded co-$t$-structures on $K^b(\operatorname{proj}\Lambda)$;
\item equivalence classes of simple-minded collections in $\Db(\Lambda)$.
\end{enumerate}
For $\Lambda = \Lambda(r,n,m)$, this bijection extends to all basic silting objects of $\Db(\Lambda)$ \cite[Theorem~6.12]{BPP2017}.
\end{theorem}

In the specific context of derived-discrete algebras, Broomhead, Pauksztello, and Ploog applied this correspondence to verify that all hearts generated by silting objects have finite length, enabling the use of abstract poset techniques.

\begin{proposition}\label{prop:algebraic-t-str}
Let $\Lambda = \Lambda(r,n,m)$. Every basic silting object in $\Db(\Lambda)$ gives rise to a bounded $t$-structure with length heart via \Cref{thm:KY}. Furthermore, because $\Db(\Lambda)$ is silting-discrete (\Cref{thm:AH}), every closed interval in the silting poset is finite.
\end{proposition}

The contractibility of the silting-pairs CW complex $\operatorname{silt}_2(\Db(\Lambda))$ is the main result of the companion paper \cite{BPP2016}.

We close this section by noting that while \cite{BPP2017,BPP2016} establish the contractibility of the silting poset and classify silting objects, they do not provide a coordinate-wise method to compute when two indecomposable objects form a partial silting pair. This motivates the hammock coordinates introduced in the subsequent section.


\section{Hammock coordinates and the pairwise partial silting criterion}
\label{sec:pairwise}

The primary goal of this section is to extract explicit integer inequalities from the Bobi\'nski--Gei\ss{}--Skowro\'nski Hom-hammock formula. By evaluating the positive suspension shifts, we translate the homological vanishing condition $\operatorname{Hom}^{>0}(Z_1, Z_2) = 0$ into a closed-form combinatorial criterion on a two-dimensional lattice. We establish this criterion first for pairs within a single component, and subsequently generalise it across distinct components.

Throughout this section, fix $\Lambda = \Lambda(r,n,m)$ with $n > r \ge 1$ and $m \ge 0$. To simplify the lattice arithmetic, we define the structural constants
\[
  a := r+m \ge 1, \qquad b := n-r \ge 1.
\]

\subsection{Suspended hammock conditions in \texorpdfstring{$(U,V)$}{(U,V)}-coordinates}\label{subsec:pairwise-susp}

For any two indecomposable objects $Z_1 = Z_{i_1, j_1}^{k_1}$ and $Z_2 = Z_{i_2, j_2}^{k_2}$ in $\operatorname{ind}(\mathcal Z)$, we define their \emph{coordinate differences} by
\[
  U := i_1 - i_2, \qquad V := j_2 - j_1.
\]
Say that a pair $\{Z_1, Z_2\}$ is a partial silting pair(=a partial silting configuration of size two) if $Z_1 \not\cong Z_2$ and $\operatorname{Hom}^{>0}(Z_1, Z_2) = \operatorname{Hom}^{>0}(Z_2, Z_1) = 0$.

We first express the unshifted degree-$0$ morphisms in terms of $(U,V)$-coordinates.

\begin{lemma}\label{lem:hammock-0}
Let $Z_1 = Z_{i_1, j_1}^k$ and $Z_2 = Z_{i_2, j_2}^k$ be objects in the same component $\mathcal Z^k$. Then the forward hammock condition of \Cref{thm:hammock}\textup{(i)} holds if and only if $U \le 0$ and $V \ge 0$ (which characterises $\operatorname{Hom}(Z_1, Z_2) \neq 0$ when $r \ge 2$, while for $r = 1$ one also has the SW branch $U \ge 1-a$ and $V \le -b-1$ from \Cref{thm:hammock}\textup{(ii)}).
\end{lemma}
\begin{proof}
By \Cref{thm:hammock}(i), $\operatorname{Hom}(Z_1, Z_2) \neq 0$ if and only if $i_2 \ge i_1$ and $j_2 \ge j_1$. This is strictly equivalent to $i_1 - i_2 \le 0$ and $j_2 - j_1 \ge 0$, which yields $U \le 0$ and $V \ge 0$.
\end{proof}

To evaluate partial silting, we must test all positive suspensions $s \ge 1$. The non-vanishing of $\operatorname{Hom}(Z_1, \Sigma^s Z_2)$ reduces to solving specific systems of inequalities parameterised by the shift $s$. To avert ambiguity at the component boundaries, we work directly in the unrolled coordinates of \eqref{eq:sigma-action}.

\begin{lemma}\label{lem:hammock-susp}
Let $Z_1, Z_2 \in \mathcal Z^k$. For any integer $s \ge 1$, $\operatorname{Hom}(Z_1, \Sigma^s Z_2) \neq 0$ if and only if one of the following two conditions holds:
\begin{enumerate}[label=\textup{(\roman*)}]
  \item \emph{NE branch:} $s = q r$ for some integer $q \ge 1$, satisfying $q a \ge U$ and $q b \le V$;
  \item \emph{SW branch:} $s = p r + 1$ for some integer $p \ge 0$, satisfying $p a \le U - 1$ and $p b \ge V + 1$.
\end{enumerate}
\end{lemma}
\begin{proof}
By \eqref{eq:sigma-action}, the suspension $\Sigma^s$ maps $\mathcal Z^k$ to $\mathcal Z^{(k+s)\bmod r}$. By \Cref{thm:hammock}, non-zero morphisms from $Z_1 \in \mathcal Z^k$ only exist into $\mathcal Z^k$ and $\mathcal Z^{k+1}$.

For the NE branch, $\Sigma^s Z_2 \in \mathcal Z^k$ if and only if $s = qr$ for some $q \ge 1$. Since $\Sigma^r$ acts on each $\mathcal Z^k$ by $(i,j)\mapsto(i+a,j-b)$, the exact coordinates are $\Sigma^{qr} Z_2 = Z_{i_2 + qa,\, j_2 - qb}^k$. By \Cref{thm:hammock}(i), non-vanishing requires $i_2 + qa \ge i_1$ and $j_2 - qb \ge j_1$, which reorganises precisely to $q a \ge U$ and $q b \le V$.

For the SW branch, $\Sigma^s Z_2 \in \mathcal Z^{k+1}$ (with cyclic identification if $k=r-1$) requires $s - 1 = pr$ for some $p \ge 0$. Writing $\Sigma^{pr+1} Z_2 = \Sigma(Z_{i_2 + pa,\, j_2 - pb}^k)$, \Cref{thm:hammock}(ii) dictates $i_2 + pa \le i_1 - 1$ and $j_2 - pb \le j_1 - 1$. This reorganises exactly to $p a \le U - 1$ and $p b \ge V + 1$.
\end{proof}

\subsection{The intra-component pairwise criterion for general \texorpdfstring{$r \ge 1$}{r >= 1}}\label{subsec:pairwise-intra}

We now solve the systems from \Cref{lem:hammock-susp} globally over all $s \ge 1$ to obtain a closed-form criterion for any pair of objects within the same component.

\begin{theorem}[Intra-component pairwise partial silting criterion]\label{thm:pairwise}
Let $\Lambda = \Lambda(r,n,m)$ with $n > r \ge 1$ and $m \ge 0$, setting $a = r+m$ and $b = n-r$. Let $Z_1, Z_2$ be distinct indecomposable objects in $\mathcal Z^k$, oriented so that $V = j_2 - j_1 \ge 0$, and set $U = i_1 - i_2$. Then $\{Z_1, Z_2\}$ is a partial silting pair if and only if one of the following systems holds:
\begin{enumerate}[label=\textup{(\Alph*)}]
  \item $V = 0$ and $1 \le |U| \le a$;
  \item $V \ge 1$, $U \ge 0$, and neither of the following conditions has an integer solution:
  \begin{align}
    \text{\textup{(I)}} &\quad \text{there exists } q \ge 1 \text{ such that } q a \ge U \text{ and } q b \le V, \label{eq:sys-I} \\
    \text{\textup{(II)}} &\quad \text{there exists } p \ge 1 \text{ such that } p a \le U - 1 \text{ and } p b \ge V + 1. \label{eq:sys-II}
  \end{align}
\end{enumerate}
Equivalently, in closed form, $\{Z_1, Z_2\}$ is a partial silting pair if and only if one of the following three conditions holds:
\begin{enumerate}[label=\textup{(\arabic*)}]
  \item[\textup{(A)}] $V = 0$ and $1 \le |U| \le a$;
  \item[\textup{(B1)}] $1 \le V \le b - 1$ and $0 \le U \le a$;
  \item[\textup{(B2)}] $V \ge b$ and $\lfloor V/b \rfloor a + 1 \le U \le (\lfloor V/b \rfloor + 1)a$, equivalently $\lfloor (U-1)/a \rfloor = \lfloor V/b \rfloor$.
\end{enumerate}
\end{theorem}

\begin{proof}
By definition, $\{Z_1, Z_2\}$ is a partial silting pair if and only if both $\operatorname{Hom}^{>0}(Z_1, Z_2) = 0$ and $\operatorname{Hom}^{>0}(Z_2, Z_1) = 0$. We apply the conditions of \Cref{lem:hammock-susp} to both $(Z_1, Z_2)$ with coordinate differences $(U, V)$ and $(Z_2, Z_1)$ with differences $(-U, -V)$.

\emph{Case $V = 0$.} Since $Z_1 \neq Z_2$, $U \neq 0$; without loss assume $U \ge 1$. For $(Z_1, Z_2)$, the NE branch requires $q b \le 0$, which is impossible since $q \ge 1$ and $b \ge 1$. The SW branch requires $p a \le U - 1$ and $p b \ge 1$. Since $b \ge 1$, the second inequality forces $p \ge 1$. The first inequality $p a \le U - 1$ has a solution $p \ge 1$ if and only if $a \le U - 1$, i.e., $U \ge a + 1$. For the reverse pair $(Z_2, Z_1)$ where $(u, v) = (-U, 0)$, the NE branch requires $q b \le 0$ (impossible) and the SW branch requires $p a \le -U - 1 < 0$ (impossible for $p \ge 0$). Thus the pair is partial silting if and only if $U \le a$, giving $1 \le |U| \le a$. This proves (A).

\emph{Case $V \ge 1$.} For the reverse pair $(Z_2, Z_1)$ with $(-U, -V)$, the NE branch $q b \le -V \le -1$ has no solution $q \ge 1$. The SW branch requires $p a \le -U - 1$ and $p b \ge -V + 1$ for $p \ge 0$. Since $V \ge 1$, $p = 0$ satisfies $0 \ge 1 - V$; at $p=0$, the first inequality requires $0 \le -U - 1 \implies U \le -1$. If $U \ge 0$, then $-U - 1 < 0 \le p a$ for all $p \ge 0$, so the SW branch fails entirely. Hence $\operatorname{Hom}^{>0}(Z_2, Z_1) = 0$ if and only if $U \ge 0$.

Now assume $U \ge 0$ and test $(Z_1, Z_2)$ with $(U, V)$. In the SW branch, $p = 0$ requires $0 \ge V + 1 \ge 2$, impossible, so any solution must have $p \ge 1$. Thus $\operatorname{Hom}^{>0}(Z_1, Z_2) = 0$ if and only if neither \eqref{eq:sys-I} nor \eqref{eq:sys-II} has a solution. This completes the proof of the system form.

Claim: The closed form is strictly equivalent to the system form.
Let $k = \lfloor V/b \rfloor \ge 0$.

In \eqref{eq:sys-I}, $q b \le V$ with $q \ge 1$ is equivalent to $1 \le q \le k$. If $1 \le V \le b - 1$, then $k = 0$, so \eqref{eq:sys-I} has no solution. If $V \ge b$, then $k \ge 1$, and $q a \ge U$ holds for some $1 \le q \le k$ if and only if it holds for the largest possible $q$, namely $q = k$. Negating $k a \ge U$ gives $U \ge k a + 1 = \lfloor V/b \rfloor a + 1$.

In \eqref{eq:sys-II}, write $V = k b + r'$ with $0 \le r' \le b - 1$. Then $p b \ge V + 1 = k b + (r' + 1)$ holds if and only if $p \ge k + 1$. Thus $p a \le U - 1$ holds for some $p \ge k + 1$ if and only if it holds for the smallest possible $p$, namely $p = k + 1$. Negating $(k + 1)a \le U - 1$ gives $U \le (k + 1)a = (\lfloor V/b \rfloor + 1)a$.

Combining these bounds with $U \ge 0$ yields exactly (B1) when $k = 0$ and (B2) when $k \ge 1$. Hence our claim.
\end{proof}

\begin{figure}[htbp]
\centering
\begin{tikzpicture}[scale=0.65]
  \draw[step=1cm, gray!25, very thin] (-2.5,-1.5) grid (8.5,7.5);

  \draw[->, thick] (-2.5,0) -- (8.5,0) node[right] {$U$};
  \draw[->, thick] (0,-1.5) -- (0,7.5) node[above] {$V$};

  \fill[red!15] (-2.5,2) rectangle (2,7.5);
  \fill[red!15] (2,4) rectangle (4,7.5);
  \fill[red!15] (4,6) rectangle (6,7.5);

  \fill[blue!15] (3,-1.5) rectangle (8.5,1);
  \fill[blue!15] (5,1) rectangle (8.5,3);
  \fill[blue!15] (7,3) rectangle (8.5,5);

  \foreach \x/\y in {-2/0,-1/0,1/0,2/0,0/1,1/1,2/1,3/2,4/2,3/3,4/3,5/4,6/4,5/5,6/5,7/6,8/6,7/7,8/7} {
    \fill[black] (\x,\y) circle (2.5pt);
  }

  \node[red, font=\bfseries] at (0,5.5) {NE forbidden};
  \node[blue, font=\bfseries] at (6.5,0) {SW forbidden};
 \node[black, font=\bfseries] at (1.8,2.5) {Allowed dots};
 
\end{tikzpicture}
\caption{Universal diagram of the $(U,V)$-lattice for $\Lambda(r,n,m)$ illustrating the intra-component pairwise criterion for $a=2, b=2$. The black dots are exactly the allowed $(U,V)$ with $0\le V\le 7$ and $-2\le U\le 8$; the shaded regions indicate the NE and SW forbidden cones.}
\label{fig:lattice}
\end{figure}
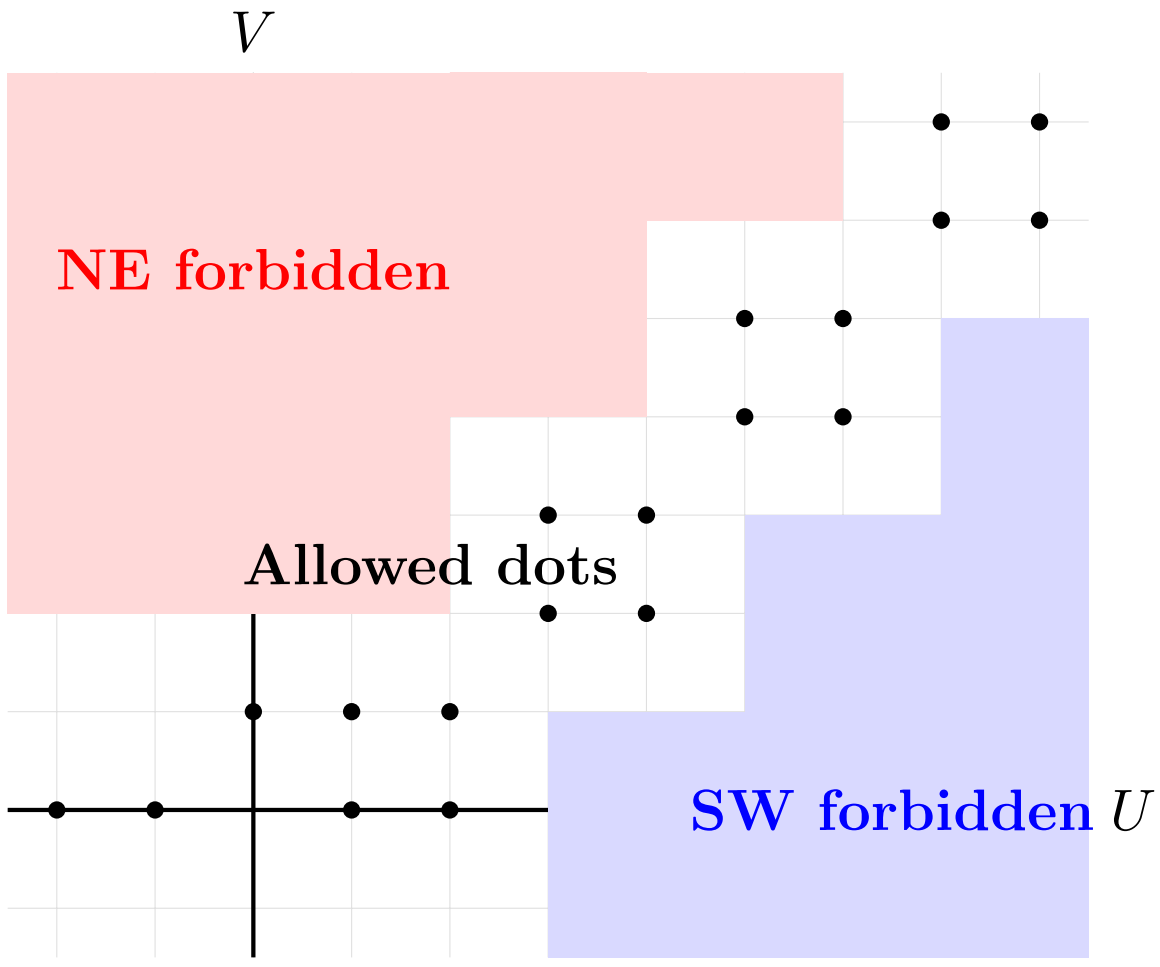

\subsection{The cross-component pairwise criterion for \texorpdfstring{$r \ge 2$}{r >= 2}}\label{subsec:pairwise-cross}

For $r \ge 2$, a basic silting object generally distributes its $\mathcal{Z}$-summands across multiple components. We extend our coordinate criterion to characterise arbitrary cross-component pairs.

\begin{theorem}[Cross-component pairwise partial silting criterion]\label{thm:cross-pairwise}
Let $\Lambda = \Lambda(r,n,m)$ with $r \ge 2$, $n>r$, and $m\ge0$. Let $Z_1 \in \mathcal Z^{k_1}$ and $Z_2 \in \mathcal Z^{k_2}$ be indecomposable objects with $0 \le k_1 < k_2 \le r-1$. Set $U = i_1 - i_2$ and $V = j_2 - j_1$. Then $\{Z_1, Z_2\}$ is a partial silting pair if and only if exactly one of the following holds:
\begin{enumerate}[label=\textup{(\roman*)}]
    \item $V = 0$ and $-a \le U \le -1$;
    \item $V \ge 1$, $U \ge 0$, and $(U,V)$ satisfies condition \textup{(B1)} or \textup{(B2)} of \Cref{thm:pairwise};
    \item $-V \ge 1$, $-U \ge 1$, and $(-U,-V)$ satisfies condition \textup{(B1)} or \textup{(B2)} of \Cref{thm:pairwise}.
\end{enumerate}
\end{theorem}

\begin{proof}
Since $0 \le k_1 < k_2 \le r-1$, any positive suspension $\Sigma^s$ ($s \ge 1$) mapping $Z_2 \in \mathcal Z^{k_2}$ into $\mathcal Z^{k_1}$ or $\mathcal Z^{(k_1+1)\bmod r}$ is of the form $s = qr + (k_1 - k_2)$ with $q \ge 1$ or $s = pr + (k_1 - k_2) + 1$ with $p \ge 1$, crossing the component boundary $\mathcal Z^{r-1} \to \mathcal Z^0$ at least once. Conversely, any suspension $\Sigma^s$ ($s \ge 1$) mapping $Z_1 \in \mathcal Z^{k_1}$ into $\mathcal Z^{k_2}$ or $\mathcal Z^{(k_2+1)\bmod r}$ is of the form $s = q'r + (k_2 - k_1)$ with $q' \ge 0$ or $s = p'r + (k_2 - k_1) + 1$ with $p' \ge 0$. Unrolling the four possible branches via \Cref{thm:hammock}, the non-vanishing conditions become
\[
\begin{array}{ll}
\text{(I)} & \exists q\ge1:\ qa\ge U,\ qb\le V, \\[2pt]
\text{(II)} & \exists p\ge1:\ pa\le U-1,\ pb\ge V+1, \\[2pt]
\text{(III)} & \exists q'\ge0:\ q'a\ge -U,\ q'b\le -V, \\[2pt]
\text{(IV)} & \exists p'\ge0:\ p'a\le -U-1,\ p'b\ge -V+1.
\end{array}
\]
Here \textup{(I)} and \textup{(II)} govern morphisms $Z_1\to \Sigma^s Z_2$ (NE and SW branches), while \textup{(III)} and \textup{(IV)} govern morphisms $Z_2\to \Sigma^s Z_1$. The pair is partial silting if and only if none of \textup{(I)}--\textup{(IV)} has a solution.

\emph{Case $V=0$.} Condition \textup{(I)} is impossible since $qb\le 0$ whereas $q \ge 1$ and $b \ge 1$. Condition \textup{(II)} has a solution if and only if $U\ge a+1$: with $p\ge1$, $pb\ge1$ is automatic, and $pa\le U-1$ has a solution iff $a\le U-1$. Condition \textup{(III)} has a solution if and only if $U\ge0$: the choice $q'=0$ gives $0\ge -U$ and $0\le 0$, so the pair fails. Condition \textup{(IV)} has a solution if and only if $U\le -a-1$: $p'=0$ fails the inequality $0\ge1$, so $p'\ge1$, and then $p'a\le -U-1$ has a solution iff $a\le -U-1$. Hence none of \textup{(I)}--\textup{(IV)} has a solution exactly when $U\le a$, $U\le-1$, and $U\ge -a$. Combining these yields $-a\le U\le-1$. This proves (i).

\emph{Case $V\ge1$.} Condition \textup{(III)} has no solution, since $q'b\le-V\le-1$ fails for every $q'\ge0$. Condition \textup{(IV)} has a solution if and only if $U\le-1$: for $U\le-1$, $p'=0$ works, and for $U\ge0$ the first inequality $p'a\le-U-1<0$ fails for every $p'\ge0$. Hence if $U\le-1$ the pair fails. Assume therefore $U\ge0$. Then \textup{(IV)} is impossible, and none of \textup{(I)}--\textup{(IV)} has a solution exactly when neither \textup{(I)} nor \textup{(II)} has a solution. Since $V\ge1$ and $U\ge0$, this is precisely the system form of \Cref{thm:pairwise}, equivalently \textup{(B1)} or \textup{(B2)}. This proves (ii).

\emph{Case $V\le-1$.} If $U\ge0$, then \textup{(III)} has the solution $q'=0$ (since $0\ge -U$ and $0\le -V$), so the pair fails. If $U\le-1$, then \textup{(I)} is impossible (since $qb\le V\le-1$ fails for $q\ge1$), and \textup{(II)} is impossible (since $pa\le U-1\le-2$ fails for $p\ge1$). Set $\widetilde U=-U\ge1$ and $\widetilde V=-V\ge1$. In \textup{(III)}, $q'=0$ fails $0 \ge \widetilde U \ge 1$, and in \textup{(IV)}, $p'=0$ fails $0 \ge \widetilde V + 1 \ge 2$; thus $q', p' \ge 1$, and \textup{(III)} and \textup{(IV)} become exactly \textup{(I)} and \textup{(II)} for $(\widetilde U,\widetilde V)$. Hence the pair is partial silting if and only if $(\widetilde U,\widetilde V)$ satisfies \textup{(B1)} or \textup{(B2)} of \Cref{thm:pairwise}. This proves (iii).
\end{proof}

\subsection{Specialisations to \texorpdfstring{$m=0$}{m=0} and small parameters}\label{subsec:pairwise-spec}

When the tail length $m=0$, the criterion collapses into a purely linear threshold without nested floor divisions.

\begin{corollary}\label{cor:case-m0}
Let $\Lambda = \Lambda(r,n,0)$, setting $a=r$ and $b=n-r$. For two distinct objects $Z_1, Z_2$ in the same component oriented with $V \ge 0$, $\{Z_1, Z_2\}$ is a partial silting pair if and only if one of the following holds:
\begin{enumerate}[label=\textup{(\Alph*)}]
  \item $V = 0$ and $1 \le |U| \le r$;
  \item $1 \le V \le n-r-1$ and $0 \le U \le r$;
  \item $V \ge n-r$ and $\lfloor (U-1)/r \rfloor = \lfloor V/(n-r) \rfloor$.
\end{enumerate}
\end{corollary}

It is readily verified that small examples yield clean numeric boundaries. We document the explicit specialisation for $\Lambda(1,3,1)$, which informs the sharp configurations of \S 6.

\begin{lemma}\label{lem:131-spec}
Let $\Lambda = \Lambda(1,3,1)$, yielding $a=2$ and $b=2$. A pair $\{Z_1, Z_2\}$ in $\mathcal Z^0$ with $V \ge 1$ is a partial silting pair if and only if one of the following holds:
\begin{enumerate}[label=\textup{(\roman*)}]
  \item $V = 1$ and $U \in \{0, 1, 2\}$;
  \item $V \in \{2k, 2k+1\}$ for some $k \ge 1$, and $U \in \{2k+1, 2k+2\}$.
\end{enumerate}
\end{lemma}
\begin{proof}
By \Cref{thm:pairwise}(B1), for $V=1 \le b-1$, we require $0 \le U \le 2$, giving (i). By \Cref{thm:pairwise}(B2), for $V \ge 2$, we require $\lfloor V/2 \rfloor \cdot 2 + 1 \le U \le (\lfloor V/2 \rfloor + 1) \cdot 2$. Letting $k = \lfloor V/2 \rfloor \ge 1$, we obtain $2k+1 \le U \le 2k+2$, giving (ii).
\end{proof}

\Cref{ex:running} (continued). We tabulate the criterion for the first two non-trivial cases to concretely illustrate the boundaries.

\begin{table}[htbp]
\centering
\small
\begin{tabular}{c|c|c|c|c|c}
\multicolumn{6}{c}{$\Lambda(1,2,0)$: $a=1, b=1 \implies V=0\Rightarrow |U|=1$; $V\ge 1\Rightarrow U=V+1$} \\
\hline
$(U,V)$ & $(0,0)$ & $(1,0)$ & $(2,1)$ & $(3,2)$ & $(1,1)$ \\
Status & $U=0$ excluded & partial & partial & partial & $U\neq V+1$
\end{tabular}
\end{table}

\begin{table}[htbp]
\centering
\small
\begin{tabular}{c|c|c|c|c|c}
\multicolumn{6}{c}{$\Lambda(1,3,0)$: $a=1, b=2 \implies V=1\Rightarrow U\in\{0,1\}$; $V\ge 2\Rightarrow U=\lfloor V/2\rfloor+1$} \\
\hline
$(U,V)$ & $(1,1)$ & $(2,2)$ & $(2,3)$ & $(3,4)$ & $(3,3)$ \\
Status & partial & partial & partial & partial & excluded
\end{tabular}
\end{table} \hfill $\triangle$


\section{The \texorpdfstring{$L$}{L}-invariant, defect arithmetic, and the torsion invariant}\label{sec:L-defect}

The primary objective of this section is to quantify the combinatorial criterion developed in \S 3. We introduce a coordinate valuation, the $L$-invariant, to associate an integer to each $\mathcal Z$-object. The difference of $L$-invariants yields the coordinate defect $\partial$, which is additive along sequences of objects. We also define the torsion invariant $g$, which dictates the minimal positive defect and later bounds the cyclic symmetry of the configuration lattice.

Throughout this section, fix $\Lambda = \Lambda(r,n,m)$ with $n > r \ge 1$ and $m \ge 0$, and set $a = r+m$ and $b = n-r$.

\subsection{The \texorpdfstring{$L$}{L}-invariant and its suspension transformation}\label{subsec:L-inv}

We assign a single integer to each indecomposable object in the $\mathcal Z$-components by taking a weighted sum of its hammock coordinates.

\begin{definition}\label{def:L}
Say that the $L$-invariant of an indecomposable object $Z_{i,j}^k \in \operatorname{ind}(\mathcal Z^k)$ is the integer
\[
  L(Z_{i,j}^k) := b\,i + a\,j = (n-r)i + (r+m)j.
\]
\end{definition}

The $L$-invariant possesses a remarkable invariance under the suspension functor $\Sigma$.

\begin{lemma}[Properties of the $L$-invariant]\label{lem:L-props}
For any indecomposable object $Z \in \operatorname{ind}(\mathcal Z)$ and any integer $s \in \mathbb Z$, we have $L(\Sigma^s Z) = L(Z)$. Consequently, $L$ is a well-defined invariant on the $\Sigma$-orbits of $\mathcal Z$-objects.
\end{lemma}

\begin{proof}
Let $Z = Z_{i,j}^k \in \mathcal Z^k$. We first evaluate the action of $\Sigma$ in two cases according to \eqref{eq:sigma-action}.

If $0 \le k < r-1$, then $\Sigma Z_{i,j}^k = Z_{i,j}^{k+1}$. Since the coordinates $(i,j)$ are strictly preserved, $L(\Sigma Z) = b\,i + a\,j = L(Z)$.

If $k = r-1$, then $\Sigma Z_{i,j}^{r-1} = Z_{i+a,\, j-b}^0$. We compute the $L$-invariant of this shifted object:
\[
  L(\Sigma Z_{i,j}^{r-1}) = b(i+a) + a(j-b) = b\,i + b\,a + a\,j - a\,b = b\,i + a\,j = L(Z_{i,j}^{r-1}).
\]
These two cases cover all possible $k$; for $r=1$ the first case is empty and the second applies. Hence $L(\Sigma Z)=L(Z)$.

Since $\Sigma$ is an autoequivalence, applying the same identity to $\Sigma^{-1}Z$ gives
\[
  L(Z) = L(\Sigma(\Sigma^{-1}Z)) = L(\Sigma^{-1}Z).
\]
Thus $L(\Sigma^s Z)=L(Z)$ for all $s\in\mathbb Z$ by induction. Hence $L$ is constant across any $\Sigma$-orbit.
\end{proof}

\subsection{The defect \texorpdfstring{$\partial$}{partial} and its additivity}\label{subsec:defect}

We can measure the homological ``distance'' between two objects by taking the difference of their $L$-invariants.

\begin{definition}\label{def:defect}
Given any two indecomposable objects $Z_1, Z_2 \in \operatorname{ind}(\mathcal Z)$, say that the coordinate defect of the ordered pair $(Z_1, Z_2)$ is the integer
\[
  \partial(Z_1, Z_2) := L(Z_1) - L(Z_2).
\]
\end{definition}

\begin{lemma}\label{lem:defect-L}
Let $Z_1 = Z_{i_1, j_1}^{k_1}$ and $Z_2 = Z_{i_2, j_2}^{k_2}$ be objects in $\operatorname{ind}(\mathcal Z)$. If we set $U = i_1 - i_2$ and $V = j_2 - j_1$, then
\[
  \partial(Z_1, Z_2) = b\,U - a\,V.
\]
\end{lemma}

\begin{proof}
By \Cref{def:L}, we expand the difference:
\[
  \partial(Z_1, Z_2) = (b\,i_1 + a\,j_1) - (b\,i_2 + a\,j_2) = b(i_1 - i_2) - a(j_2 - j_1) = b\,U - a\,V.
\]
This completes the proof.
\end{proof}

The defect is computationally well-behaved and inherits strict structural properties from the $L$-invariant.

\begin{lemma}[Properties of the defect]\label{lem:defect-props}
For any objects $Z_1, Z_2 \in \operatorname{ind}(\mathcal Z)$ and any integer $s \in \mathbb Z$, the following hold:
\begin{enumerate}[label=\textup{(\roman*)}]
    \item \emph{Skew-symmetry:} $\partial(Z_2, Z_1) = -\partial(Z_1, Z_2)$;
    \item \emph{$\Sigma$-invariance:} $\partial(\Sigma^s Z_1, \Sigma^s Z_2) = \partial(Z_1, Z_2)$;
    \item \emph{Integrality:} $\partial(Z_1, Z_2) \in \mathbb Z$.
\end{enumerate}
\end{lemma}

\begin{proof}
Skew-symmetry and integrality follow directly from \Cref{def:defect} and integer arithmetic. For $\Sigma$-invariance, \Cref{lem:L-props} gives $L(\Sigma^s Z_p) = L(Z_p)$ for $p \in \{1, 2\}$ and all $s\in\mathbb Z$, hence
\[
  \partial(\Sigma^s Z_1, \Sigma^s Z_2) = L(\Sigma^s Z_1) - L(\Sigma^s Z_2) = L(Z_1) - L(Z_2) = \partial(Z_1, Z_2).
\]
This completes the proof.
\end{proof}

The most crucial property of the defect is its additivity. Unlike the complex non-linear intervals of the partial silting criteria in \S 3, the defect accumulates linearly.

\begin{corollary}[Additivity along chains]\label{cor:defect-chain}
Suppose $Z_1, Z_2, \ldots, Z_k$ is a sequence of objects in $\operatorname{ind}(\mathcal Z)$. Then the defect of the outer pair is exactly the sum of the defects of the consecutive steps:
\[
  \partial(Z_1, Z_k) = \sum_{p=1}^{k-1} \partial(Z_p, Z_{p+1}).
\]
\end{corollary}

\begin{proof}
By \Cref{def:defect}, $\partial(Z_p, Z_{p+1}) = L(Z_p) - L(Z_{p+1})$. The summation is a telescoping series:
\[
  \sum_{p=1}^{k-1} (L(Z_p) - L(Z_{p+1})) = L(Z_1) - L(Z_k) = \partial(Z_1, Z_k).
\]
The conclusion follows.
\end{proof}

In view of the above corollary, analyzing the defect of consecutive pairs in a partial silting configuration controls the global defect. For $m=0$, we classify these consecutive steps into three distinct types based on the criteria in \Cref{cor:case-m0}.

\begin{definition}[Step types for $m=0$]\label{def:step-types}
Let $\Lambda = \Lambda(r,n,0)$, so $a=r$ and $b=n-r$. Given a partial silting pair $\{Z_1, Z_2\}$ oriented such that $V = j_2 - j_1 \ge 0$, say that the step from $Z_1$ to $Z_2$ is of:
\begin{itemize}
    \item \textbf{Type A:} if $V = 0$. By \Cref{cor:case-m0}(A), $1 \le |U| \le r$, so $|\partial(Z_1, Z_2)| = b|U| \ge b \ge 1$.
    \item \textbf{Type B+:} if $V \ge 1$ and $U \ge 1$.
    \item \textbf{Type B0:} if $V \ge 1$ and $U = 0$. By \Cref{cor:case-m0}\textup{(B1)}, this occurs only when $1 \le V \le b-1$. Crucially, $\partial(Z_1, Z_2) = -aV < 0$.
\end{itemize}
By \Cref{cor:case-m0}, every partial silting pair with $V\ge0$ falls into exactly one of these three types.
\end{definition}

\subsection{The torsion invariant and minimal positive defect}
\label{subsec:torsion}

We define a final structural constant derived from the parameters $a$ and $b$, which restricts the possible values the defect $\partial$ can attain.

\begin{definition}\label{def:torsion}
Say that the torsion invariant $g$ of the derived-discrete algebra $\Lambda(r,n,m)$ is
\[
  g := \gcd(r+m, n-r) = \gcd(a, b).
\]
\end{definition}

\begin{proposition}[Minimal positive defect]\label{prop:min-defect}
The set of all possible defects
\[
  \{\partial(Z_1, Z_2) \mid Z_1, Z_2 \in \operatorname{ind}(\mathcal Z)\}
\]
is exactly the principal ideal $g \mathbb Z$. In particular, the minimal positive defect between any two $\mathcal Z$-objects is exactly $g$.
\end{proposition}

\begin{proof}
Let
\[
  I = \{\partial(Z_1, Z_2) \mid Z_1, Z_2 \in \operatorname{ind}(\mathcal Z)\}.
\]

Claim: $I = \{b\,U - a\,V \mid U, V \in \mathbb Z\}$.
Given any $U, V \in \mathbb Z$, choose $Z_1 = Z_{U, 0}^0$ and $Z_2 = Z_{0, V}^0$ in $\mathcal Z^0$; this is possible because the coordinates $(i,j)$ range over all integers in the $\mathcal Z$-components. Then
\[
  i_1 - i_2 = U - 0 = U, \qquad j_2 - j_1 = V - 0 = V.
\]
By \Cref{lem:defect-L}, $\partial(Z_1, Z_2) = b\,U - a\,V$. Conversely, every pair of objects yields some $U, V \in \mathbb Z$. Hence our claim.

By B\'ezout's identity, the set of all integer linear combinations
\[
  \{b\,U - a\,V \mid U, V \in \mathbb Z\}
\]
is exactly the ideal generated by the greatest common divisor of $a$ and $b$. Thus
\[
  I = \gcd(a,b)\mathbb Z = g \mathbb Z.
\]
It is readily verified that the smallest strictly positive element of this ideal is $g$.
\end{proof}


\section{The universal upper bound: combinatorial staircase and Bongartz completion}
\label{sec:bound}

The main goal of this section is to establish the universal upper bound $k \le n+m$ on the size of any partial silting configuration of $\mathcal Z$-objects in $\Db(\Lambda(r,n,m))$. We first show that pairwise partial silting guarantees that the direct sum is a presilting object. In view of the highly complex nature of the derived category, we provide a purely combinatorial proof of the upper bound within a single component using a monotonic $\Sigma^r$-invariant staircase. We close this section by extending the bound globally across all components via Bongartz completion, illustrating the result with a complete classification of silting configurations for $\Lambda(1,2,0)$.

Throughout this section, we fix $\Lambda = \Lambda(r,n,m)$ with $n > r \ge 1$ and $m \ge 0$, and set $a = r+m$ and $b = n-r$.

\subsection{Self-extension vanishing and presilting direct sums}
\label{subsec:bound-presilting}

Let $S = \{X_1, \ldots, X_k\}$ be a set of pairwise non-isomorphic indecomposable $\mathcal Z$-objects in $\Db(\Lambda)$. Recall that $S$ is a partial silting configuration if $\operatorname{Hom}^{>0}(X_p, X_q) = 0$ for all $p \neq q$. To treat the direct sum $M := \bigoplus_{p=1}^k X_p$ as a presilting object(=an object with vanishing positive self-extensions), we must guarantee that the self-extensions of each individual object $X_p$ vanish.

\begin{lemma}\label{lem:self-hom-vanishing}
Let $X_p = Z_{i,j}^k \in \operatorname{ind}(\mathcal Z)$. For any integer $s \ge 1$, we have
\[
  \operatorname{Hom}(X_p, \Sigma^s X_p) = 0.
\]
\end{lemma}

\begin{proof}
By \eqref{eq:sigma-action}, the suspension $\Sigma^s$ maps $\mathcal Z^k$ to $\mathcal Z^{(k+s)\bmod r}$. By \Cref{thm:hammock}, non-zero morphisms originating from $\mathcal Z^k$ strictly land in $\mathcal Z^k$ or $\mathcal Z^{(k+1)\bmod r}$.

Claim: If $(k+s) \bmod r \notin \{k, (k+1) \bmod r\}$, then $\operatorname{Hom}(X_p, \Sigma^s X_p) = 0$.
This follows immediately from \Cref{thm:hammock}(iii). Hence our claim.

In view of the claim, it remains to evaluate the two adjacent components.

If $(k+s) \bmod r = k$, we write $s = q r$ with $q \ge 1$. By the unrolled suspension formula,
\[
  \Sigma^{qr} X_p = Z_{i+qa,\, j-qb}^k.
\]
By \Cref{thm:hammock}(i), non-vanishing requires $i+qa \ge i$ and $j-qb \ge j$. The latter inequality implies $qb \le 0$. Since $q \ge 1$ and $b = n-r \ge 1$, this is impossible.

If $(k+s) \bmod r = (k+1) \bmod r$, we write $s = p r + 1$ with $p \ge 0$. By \Cref{thm:hammock}(ii), non-vanishing requires
\[
  \Sigma^{pr} X_p = Z_{i+pa,\, j-pb}^k
\]
to satisfy $i+pa \le i-1$ and $j-pb \le j-1$. The first inequality gives $pa \le -1$. Since $p \ge 0$ and $a = r+m \ge 1$, this is impossible. For $r=1$, both branches land in $\mathcal{Z}^0$ and must be excluded simultaneously; since $qb \le 0$ ($q \ge 1$) and $pa \le -1$ ($p \ge 0$) both fail as shown above, $\operatorname{Hom}(X_p, \Sigma^s X_p) = 0$ holds unconditionally for all $r \ge 1$. Thus our assumption that $\operatorname{Hom}(X_p, \Sigma^s X_p) \neq 0$ is wrong, and the conclusion follows.
\end{proof}

\begin{lemma}\label{lem:presilting}
If $S = \{X_1, \ldots, X_k\}$ is a partial silting configuration of pairwise non-isomorphic $\mathcal Z$-objects, then the direct sum $M := \bigoplus_{p=1}^k X_p$ is a basic presilting object of $\Db(\Lambda)$.
\end{lemma}

\begin{proof}
Since the $X_p$ are pairwise non-isomorphic indecomposables, $M$ is basic. By additivity of the Hom-functor,
\[
  \operatorname{Hom}^{>0}(M, M)
  \cong \bigoplus_{p,q} \operatorname{Hom}^{>0}(X_p, X_q).
\]
For $p \neq q$, the Hom-spaces vanish by the definition of a partial silting configuration. For $p = q$, the spaces vanish by \Cref{lem:self-hom-vanishing}. Hence $\operatorname{Hom}^{>0}(M, M) = 0$.
\end{proof}

\subsection{Unconditional combinatorial bound and staircase classification within a component}
\label{subsec:bound-combinatorial}

While \Cref{thm:bound} below establishes the upper bound $k \le n+m$ globally across all $\mathcal Z$-components via Bongartz completion, we provide a direct combinatorial proof that any partial silting configuration contained entirely within a single component $\mathcal Z^k$ has at most $n+m$ elements. This approach yields a structural parameterization of maximal partial silting configurations in terms of monotonic staircase subsets.

For each indecomposable object $Z_p = Z_{i_p, j_p}^k \in \operatorname{ind}(\mathcal Z^k)$, we write the Euclidean division of $j_p$ by $b$ as
\[
  j_p = q_p b + r_p, \qquad
  q_p = \left\lfloor \frac{j_p}{b} \right\rfloor \in \mathbb Z,
  \quad r_p \in \{0, 1, \ldots, b-1\},
\]
and define the translated $i$-coordinate by
\[
  I_p := i_p + q_p a = i_p + \left\lfloor \frac{j_p}{b} \right\rfloor a.
\]
Since $\Sigma^r(Z_{i,j}^k) = Z_{i+a, \, j-b}^k$, both $I_p$ and $r_p$ are strictly invariant under $\Sigma^r$. Two objects in $\mathcal Z^k$ lie in the same $\Sigma^r$-orbit if and only if $(I_p, r_p) = (I_q, r_q)$.

\begin{lemma}[Monotonicity of $(I,r)$-coordinates]\label{lem:Ir-monotonicity}
Let $\{Z_p, Z_q\}$ be a partial silting pair of distinct objects in $\mathcal Z^k$.
\begin{enumerate}[label=\textup{(\roman*)}]
  \item If $r_p = r_q$, then $1 \le |I_p - I_q| \le a$. In particular, $(I_p, r_p) \neq (I_q, r_q)$.
  \item If $r_p < r_q$, then $0 \le I_p - I_q \le a$.
\end{enumerate}
\end{lemma}

\begin{proof}
Without loss we relabel $\{Z_p, Z_q\}$ as $(Z_1, Z_2)$ such that $j_1 \le j_2$. Set
\[
  U = i_1 - i_2,
  \qquad
  V = j_2 - j_1 = (q_2 - q_1)b + (r_2 - r_1) \ge 0.
\]
Then
\[
  I_1 - I_2 = U - (q_2 - q_1)a.
\]

\emph{Case $r_1 = r_2$.} Here $V = (q_2 - q_1)b$. If $q_1 = q_2$, then $V = 0$, and \Cref{thm:pairwise}(A) guarantees $1 \le |U| \le a$, yielding $1 \le |I_1 - I_2| \le a$. If $q_2 > q_1$, then $V \ge b$ with $\lfloor V/b \rfloor = q_2 - q_1$. \Cref{thm:pairwise}(B2) gives
\[
  (q_2 - q_1)a + 1 \le U \le (q_2 - q_1 + 1)a,
\]
which reorganizes exactly to $1 \le I_1 - I_2 \le a$.

\emph{Case $r_1 < r_2$.} Here $\lfloor V/b \rfloor = q_2 - q_1 \ge 0$. If $q_1 = q_2$, then $1 \le V \le b-1$, and \Cref{thm:pairwise}(B1) gives $0 \le U \le a$, so $0 \le I_1 - I_2 \le a$. If $q_2 > q_1$, then $V \ge b$, and \Cref{thm:pairwise}(B2) yields
\[
  1 \le I_1 - I_2 \le a,
\]
in particular $0 \le I_1 - I_2 \le a$.

\emph{Case $r_1 > r_2$.} Since $V > 0$ and $-(b-1) \le r_2 - r_1 \le -1$, we have $q_2 - q_1 \ge 1$ and $\lfloor V/b \rfloor = q_2 - q_1 - 1$. If $q_2 - q_1 = 1$, then $1 \le V \le b-1$, and \Cref{thm:pairwise}(B1) gives $0 \le U \le a$, whence
\[
  I_2 - I_1 = a - U \in [0, a].
\]
If $q_2 - q_1 \ge 2$, then $V \ge b$, and \Cref{thm:pairwise}(B2) gives
\[
  (q_2 - q_1 - 1)a + 1 \le U \le (q_2 - q_1)a,
\]
whence
\[
  I_2 - I_1 = (q_2 - q_1)a - U \in [0, a-1] \subseteq [0, a].
\]
Since $r_2 < r_1$, this again gives $0 \le I_{\text{smaller } r} - I_{\text{larger } r} \le a$.
\end{proof}

\begin{theorem}[Combinatorial intra-component bound]\label{thm:sharp-combinatorial}
Let $S = \{Z_1, \ldots, Z_k\} \subseteq \operatorname{ind}(\mathcal Z^\ell)$ be a partial silting configuration of pairwise non-isomorphic objects in a single component $\mathcal Z^\ell$. Then $k \le n + m$.
\end{theorem}

\begin{proof}
Order the objects $Z_{(1)}, \ldots, Z_{(k)}$ of $S$ so that
\[
  r_{(1)} \le r_{(2)} \le \cdots \le r_{(k)},
\]
breaking ties for equal residues by $I_{(t)} > I_{(t+1)}$; this is possible since $r_p = r_q \implies I_p \neq I_q$ by \Cref{lem:Ir-monotonicity}(i). For each consecutive step $1 \le t \le k-1$:
\begin{itemize}
  \item if $r_{(t)} = r_{(t+1)}$, then $r_{(t+1)} - r_{(t)} = 0$ and $I_{(t)} - I_{(t+1)} \ge 1$ by the tie-breaking order;
  \item if $r_{(t)} < r_{(t+1)}$, then $r_{(t+1)} - r_{(t)} \ge 1$ and $I_{(t)} - I_{(t+1)} \ge 0$ by \Cref{lem:Ir-monotonicity}(ii).
\end{itemize}
In either case,
\[
  \bigl(r_{(t+1)} - r_{(t)}\bigr) + \bigl(I_{(t)} - I_{(t+1)}\bigr) \ge 1.
\]
Summing from $t = 1$ to $k - 1$ telescopes to
\[
  k - 1
  \;\le\;
  \bigl(r_{(k)} - r_{(1)}\bigr) + \bigl(I_{(1)} - I_{(k)}\bigr)
  \;\le\;
  (b - 1) + a
  \;=\;
  (n - r - 1) + (r + m)
  \;=\;
  n + m - 1,
\]
where $r_{(k)} - r_{(1)} \le b - 1$ since $r_p \in \{0, \ldots, b-1\}$, and $I_{(1)} - I_{(k)} \le a$ by \Cref{lem:Ir-monotonicity} applied to the pair $\{Z_{(1)}, Z_{(k)}\}$ with $r_{(1)} \le r_{(k)}$. Thus $k \le n + m$.
\end{proof}

\begin{corollary}[Staircase classification]\label{cor:staircase-classification}
Let $S$ be a partial silting configuration in a single component $\mathcal Z^\ell$ with $|S| = n+m$. Then there exist an integer $c$ and nonnegative integers $s_0, \ldots, s_{b-1}$ with
\[
s_0 + s_1 + \cdots + s_{b-1} = a
\]
such that the $(\rho, I)$-coordinates of $S$ are exactly
\[
\bigl\{\, (\rho, c - S_\rho - t) \;\bigm|\; 0 \le \rho \le b-1,\ 0 \le t \le s_\rho \,\bigr\},
\]
where $S_\rho := s_0 + \cdots + s_{\rho-1}$ for $\rho \ge 1$ and $S_0 := 0$. Equivalently, the configuration is a staircase with $b$ columns of heights $s_\rho+1$ whose total height decrease is $a$. Conversely, every such staircase is realised by a partial silting configuration in $\mathcal Z^\ell$ by choosing all representatives with the same quotient $q$.
\end{corollary}

\begin{proof}
Set $k = n+m$. By the equality case of the telescoping argument in \Cref{thm:sharp-combinatorial}, every consecutive pair in the $(r_p, I_p)$-ordering satisfies either ($r_{(t+1)} = r_{(t)}$ and $I_{(t)} - I_{(t+1)} = 1$) or ($r_{(t+1)} = r_{(t)} + 1$ and $I_{(t)} = I_{(t+1)}$); additionally $r_{(1)} = 0$, $r_{(k)} = b-1$, and $I_{(1)} - I_{(k)} = a$.

By \Cref{lem:Ir-monotonicity}\textup{(i)}, the coordinate map $Z_p \mapsto (r_p, I_p)$ is strictly injective on $S$. Since $r_{(1)} = 0$, $r_{(k)} = b-1$, and $r_{(t+1)} - r_{(t)} \in \{0,1\}$, every residue $\rho \in \{0, \ldots, b-1\}$ appears in $S$ and hosts a contiguous vertical block of $I$-values.

Let $s_\rho \ge 0$ denote the height of column $\rho$ minus $1$; set $S_\rho = s_0 + \cdots + s_{\rho-1}$ and $c = I_{(1)}$. The consecutive-step structure forces the bottom of column $\rho$ to coincide in height with the top of column $\rho+1$, yielding top height $c - S_\rho$ in column $\rho$ by induction. The count $|S| = \sum_{\rho=0}^{b-1} (s_\rho + 1) = n+m$ yields $\sum_{\rho=0}^{b-1} s_\rho = n+m - b = a$, proving the forward direction.

For the converse, given $(c, s_0, \ldots, s_{b-1})$ with $\sum_\rho s_\rho = a$, fix a quotient $q \in \mathbb Z$ and define the objects $Z_{i_\rho(t), \, j_\rho}^\ell$ with coordinates $j_\rho = qb + \rho$ and $i_\rho(t) = c - S_\rho - t - qa$ for $0 \le \rho \le b-1$ and $0 \le t \le s_\rho$. Every pair with equal residues satisfies $V = 0$ and $1 \le |U| \le a$ (\Cref{thm:pairwise}\textup{(A)}), and every pair with distinct residues satisfies $1 \le V \le b-1$ and $0 \le U \le a$ (\Cref{thm:pairwise}\textup{(B1)}), so the resulting set is a partial silting configuration of size $n+m$.
\end{proof}

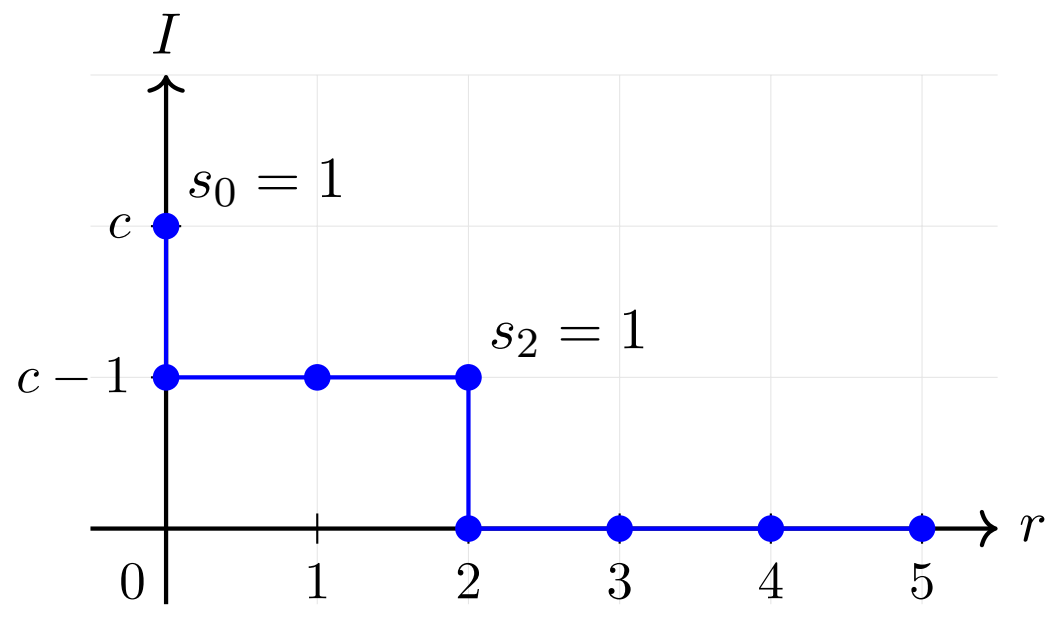
\begin{figure}[htbp]
\centering
\begin{tikzpicture}[scale=1.0]
\draw[step=1cm, gray!20, very thin] (-0.5, -0.5) grid (5.5, 3.0);

\draw[->, thick] (-0.5,0) -- (5.5,0) node[right] {$r$};
\draw[->, thick] (0,-0.5) -- (0,3.0) node[above] {$I$};

\foreach \x in {1,2,3,4,5} \draw (\x,0.1) -- (\x,-0.1) node[below] {\small $\x$};
\draw (0,0.1) -- (0,-0.1) node[below left] {\small $0$};
\draw (0.1,1) -- (-0.1,1) node[left] {\small $c-1$};
\draw (0.1,2) -- (-0.1,2) node[left] {\small $c$};

\draw[thick, blue] (0,2) -- (0,1) -- (1,1) -- (2,1) -- (2,0) -- (5,0);

\foreach \x/\y in {0/2, 0/1, 1/1, 2/1, 2/0, 3/0, 4/0, 5/0} {
\fill[blue] (\x,\y) circle (2.5pt);
}

\node[above right, font=\bfseries] at (0,2) {$s_0=1$};
\node[above right, font=\bfseries] at (2,1) {$s_2=1$};
\end{tikzpicture}
\caption{The $\Sigma^r$-invariant $(I,r)$-staircase for $a=2$ and $b=6$. The column at residue $\rho$ has height $s_\rho+1$; here $s_0=1$, $s_2=1$, and all other $s_\rho=0$, giving $a+b=8$ points.}
\label{fig:staircase}
\end{figure}

\subsection{Global upper bound via Bongartz completion}
\label{subsec:bound-bongartz}

We extend the combinatorial bound globally across all $r$ components by utilizing Bongartz completion.

\begin{theorem}[Universal upper bound]\label{thm:bound}
Let $S=\{X_1,\ldots,X_k\}$ be a partial silting configuration of pairwise non-isomorphic indecomposable $\mathcal Z$-objects in $\Db(\Lambda(r,n,m))$. Then $k \le n+m$.
\end{theorem}

\begin{proof}
By \Cref{lem:presilting}, the direct sum $M=\bigoplus_{p=1}^{k} X_p$ is a basic presilting object of $\Db(\Lambda)$. Since $\Lambda$ is derived-discrete, $\Db(\Lambda)$ is a silting-discrete triangulated category by \Cref{thm:AH}. Thus, the Bongartz completion result of \Cref{thm:AM} guarantees that the presilting object $M$ can be extended to a basic silting object $T = M\oplus M'$ of $\Db(\Lambda)$.

The non-isomorphic indecomposable summands of the basic silting object $T$ form a basis for the Grothendieck group $K_0(\Db(\Lambda))$ by \Cref{lem:silting-basic}. Since the rank of $K_0(\Db(\Lambda))$ is exactly $n+m$, it is readily verified that $|T| = n+m$. As $M$ is a direct summand of $T$, we conclude $k = |M| \le |T| = n+m$.
\end{proof}

\begin{corollary}\label{cor:bound-1n0}
For $\Lambda = \Lambda(1,n,0)$, any partial silting configuration of $\mathcal Z$-objects has at most $n$ elements.
\end{corollary}

We apply the coordinate criteria to yield a complete classification of silting configurations for the smallest algebra, $\Lambda(1,2,0)$.

\Cref{ex:running} (continued). For $\Lambda(1,2,0)$, we have $a=1, b=1$, and the maximal configuration size is $n+m = 2$. By \Cref{cor:case-m0}, a pair $\{Z_1, Z_2\}$ with $V = j_2 - j_1 \ge 0$ is a partial silting pair if and only if either ($V=0$ and $|U|=1$) or ($V \ge 1$ and $U = V+1$). It is easy to check that every maximal configuration of $\Lambda(1,2,0)$ up to $\Sigma$-equivalence falls exactly into one of three geometric step types:
\begin{itemize}
    \item \textbf{Type I-a:} $\{Z_{i,j}, Z_{i-1,j}\}$. The difference is $V=0$ and $U=1$.
    \item \textbf{Type I-b:} $\{Z_{i,j}, Z_{i+1,j}\}$. The difference is $V=0$ and $U=-1$.
    \item \textbf{Type II:} $\{Z_{i_1,j_1}, Z_{i_2,j_2}\}$ with $V \ge 1$ and $U = V+1$.
\end{itemize}
Thus, hammock coordinates recover the known silting discrete structures of the Kronecker-like gentle algebra entirely algebraically. \hfill $\triangle$


\section{Explicit sharp configurations attaining the bound \texorpdfstring{$n+m$}{n+m}}\label{sec:sharp}

The primary objective of this section is to show that the universal upper bound $n+m$ established in \S 5 is universally sharp for every derived-discrete algebra $\Lambda(r,n,m)$. The first major achievement of this section is the explicit construction of a universal L-shaped configuration that attains this bound strictly within a single $\mathcal Z$-component. Subsequently, we present specialized diagonal configurations and apply our criteria to the running examples to illustrate the diverse geometric profiles of maximal partial silting sets.

Throughout this section, fix $\Lambda = \Lambda(r,n,m)$ with $n > r \ge 1$ and $m \ge 0$, and set $a = r+m$ and $b = n-r$.

\subsection{Universal L-shaped sharp configuration for all \texorpdfstring{$\Lambda(r,n,m)$}{Lambda(r,n,m)}}\label{sec:sharp-general}

Remarkably, the bound $n+m$ is always attained within the single component $\mathcal Z^0$ by an explicit L-shaped configuration combining $a$ horizontal steps and $b$ vertical steps.

\begin{theorem}[Universal sharpness of the $n+m$ bound]\label{thm:sharp-universal}
Let $\Lambda = \Lambda(r,n,m)$ with $n > r \ge 1$ and $m \ge 0$. Define the L-shaped configuration of $a + b = n + m$ objects in $\mathcal Z^0$ by
\[
  S \;:=\; \bigl\{\, X_p = Z_{-p,\,0}^0 \;\big|\; 0 \le p \le a-1 \,\bigr\}
  \;\sqcup\;
  \bigl\{\, Y_q = Z_{-a,\,q}^0 \;\big|\; 0 \le q \le b-1 \,\bigr\}.
\]
Then $S$ is a partial silting configuration of size $n + m$ whose elements lie in pairwise distinct $\Sigma$-orbits.
\end{theorem}

\begin{proof}
Since $-a < -p$ for all $0 \le p \le a-1$, the sets $\{X_p\}_{p=0}^{a-1}$ and $\{Y_q\}_{q=0}^{b-1}$ are disjoint. Thus $|S| = a + b = n + m$.

Claim: The elements of $S$ lie in pairwise distinct $\Sigma$-orbits.
Since every object in $S$ lies in $\mathcal Z^0$, any $\Sigma$-equivalence between two objects of $S$ must arise from a shift $s = t r$ for some integer $t$, which shifts the $j$-coordinate by $-t b$ via \eqref{eq:sigma-action}. Since the $j$-coordinates of all objects in $S$ lie in $\{0, \ldots, b-1\}$, the difference $-tb$ has absolute value strictly less than $b$, forcing $t = 0$. Thus no two distinct objects in $S$ are $\Sigma$-equivalent.

It remains to verify that every pair of distinct objects in $S$ satisfies the intra-component criterion of \Cref{thm:pairwise}:
\begin{itemize}
  \item \emph{Pairs $(X_p, X_{p'})$ with $0 \le p < p' \le a-1$:} Here $j_1 = j_2 = 0$, so $V = 0$. We have $U = (-p) - (-p') = p' - p \in \{1, \ldots, a-1\} \subseteq \{1, \ldots, a\}$. By \Cref{thm:pairwise}(A), the pair is partial silting.
  \item \emph{Pairs $(Y_q, Y_{q'})$ with $0 \le q < q' \le b-1$:} Here $i_1 = i_2 = -a$ and $j_1 < j_2$, so $U = 0$ and $V = q' - q \in \{1, \ldots, b-1\}$. Since $0 \le U \le a$, the pair is partial silting by \Cref{thm:pairwise}(B1).
  \item \emph{Cross-pairs $(X_p, Y_q)$ with $0 \le p \le a-1$ and $0 \le q \le b-1$:} Since $j(X_p) = 0 \le q = j(Y_q)$, we have $V = q \in \{0, \ldots, b-1\}$ and $U = (-p) - (-a) = a - p \in \{1, \ldots, a\}$. If $q = 0$, then $V = 0$ and $1 \le U \le a$, so \Cref{thm:pairwise}(A) holds. If $1 \le q \le b-1$, then $1 \le V \le b-1$ and $0 \le U \le a$, so \Cref{thm:pairwise}(B1) holds.
\end{itemize}
Thus every pair in $S$ is partial silting.
\end{proof}

\begin{figure}[htbp]
\centering
\begin{tikzpicture}[scale=0.85, >=stealth]
\draw[step=1cm, gray!25, very thin] (-3.5,-1.0) grid (1.5, 2.8);
\draw[->, thick] (-3.5,0) -- (1.5,0) node[right] {$i$};
\draw[->, thick] (0,-1.0) -- (0,2.8) node[above] {$j$};

\foreach \x in {-3,-2,-1,1}
  \node[below, font=\scriptsize] at (\x,-0.05) {$\x$};
\foreach \y in {1,2}
  \node[left, font=\scriptsize] at (-0.05,\y) {$\y$};

\draw[thick, blue] (0,0) -- (-1,0) -- (-2,0) -- (-2,1);
\fill[blue] (0,0) circle (2.5pt) node[below right, font=\small] {$X_0$};
\fill[blue] (-1,0) circle (2.5pt) node[above, font=\small, yshift=2pt] {$X_1$};
\fill[blue] (-2,0) circle (2.5pt) node[above left, font=\small] {$Y_0$};
\fill[blue] (-2,1) circle (2.5pt) node[left, font=\small, xshift=-2pt] {$Y_1$};

\node[blue, font=\bfseries] at (-1.0, 3.4) {$\Lambda(1,3,1)$: $a=2,\,b=2$};
\end{tikzpicture}
\hspace{1.5cm}
\begin{tikzpicture}[scale=0.85, >=stealth]
\draw[step=1cm, gray!25, very thin] (-3.5,-1.0) grid (1.5, 2.8);
\draw[->, thick] (-3.5,0) -- (1.5,0) node[right] {$i$};
\draw[->, thick] (0,-1.0) -- (0,2.8) node[above] {$j$};

\foreach \x in {-3,-2,-1,1}
  \node[below, font=\scriptsize] at (\x,-0.05) {$\x$};
\foreach \y in {1,2}
  \node[left, font=\scriptsize] at (-0.05,\y) {$\y$};

\draw[thick, red, dashed] (0,0) -- (-1,1) -- (-2,2);
\fill[red] (0,0) circle (2.5pt) node[below right, font=\small] {$Z_0$};
\fill[red] (-1,1) circle (2.5pt) node[above right, font=\small] {$Z_1$};
\fill[red] (-2,2) circle (2.5pt) node[above right, font=\small] {$Z_2$};

\node[red, font=\bfseries] at (-1.0, 3.4) {$\Lambda(1,3,0)$: $a=1,\,b=2$};
\end{tikzpicture}
\caption{Explicit hammock coordinates $(i,j)$ of two maximal partial silting configurations. Left: the L-shape of \Cref{thm:sharp-universal} in $\mathcal Z^0$ for $\Lambda(1,3,1)$, of size $a+b=4=n+m$. Right: the diagonal chain of \Cref{thm:diag} in $\mathcal Z^0$ for $\Lambda(1,3,0)$, of size $n=3=n+m$. Both panels use the standard coordinates of the respective $\mathcal{Z}^0$-component; the two algebras are distinct, so no object is shared between the panels.}
\label{fig:sharp-configs}
\end{figure}
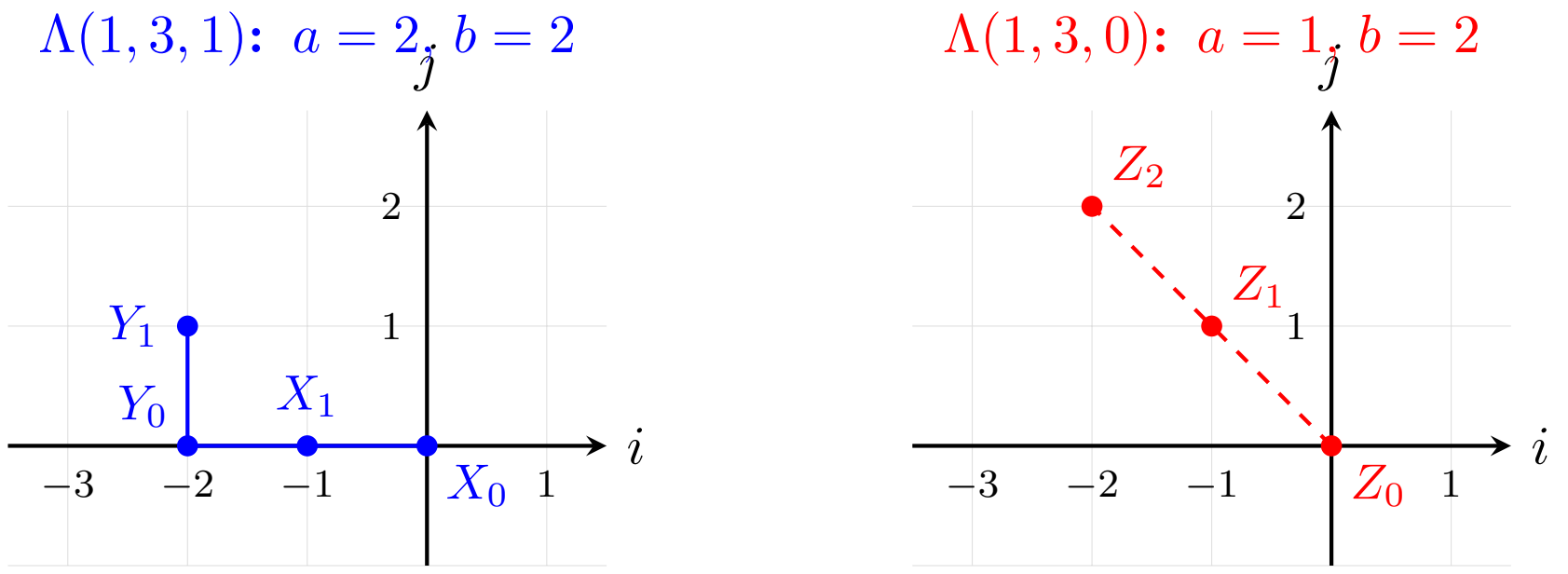
\subsection{The diagonal chain for \texorpdfstring{$\Lambda(1,n,0)$}{Lambda(1,n,0)}}\label{sec:sharp-diag}

While the L-shaped configuration is universal, the strict diagonal chain provides an alternative monotonic coordinate profile that is of independent combinatorial interest for the tail-free case $m=0$.

\begin{theorem}\label{thm:diag}
Let $\Lambda = \Lambda(1,n,0)$ with $n \ge 2$. The diagonal chain defined by
\[
  S := \bigl\{\, Z_s = Z_{-s,\, s(n-2)}^0 \;\big|\; 0 \le s \le n-1 \,\bigr\}
\]
is a partial silting configuration of size $n$, whose elements lie in pairwise distinct $\Sigma$-orbits.
\end{theorem}

\begin{proof}
Claim: The objects in $S$ lie in pairwise distinct $\Sigma$-orbits.
Suppose, for contradiction, that $Z_p$ and $Z_q$ are $\Sigma$-equivalent for some $0 \le p < q \le n-1$. Since $r=1$, the suspension $\Sigma$ shifts the coordinates by $(i,j) \mapsto (i+1, j-(n-1))$. If $\Sigma^t Z_q = Z_p$, then the $i$-coordinates dictate $(-q) + t = -p$, which forces $t = q-p$. The corresponding shift in the $j$-coordinate dictates $q(n-2) - t(n-1) = p(n-2)$. Substituting $t = q-p$ yields
\[
  q(n-2) - (q-p)(n-1) = p(n-2),
\]
equivalently
\[
  (q-p)(n-2) = (q-p)(n-1).
\]
Since $q-p > 0$, this implies $n-2 = n-1$, which is impossible. Hence our assumption is wrong, and the conclusion follows.

We now verify the partial silting conditions of \Cref{cor:case-m0} for any pair $Z_p, Z_q \in S$ with $0 \le p < q \le n-1$. The coordinate differences are
\[
  U = (-p) - (-q) = q - p,
  \qquad
  V = q(n-2) - p(n-2) = (q-p)(n-2).
\]
Note that $U = q-p \ge 1$.

If $n = 2$, then $V = (q-p)(0) = 0$. \Cref{cor:case-m0}(A) requires $1 \le |U| \le 1$, which is satisfied since $U = q-p = 1$ (as $p < q \le 1$).

If $n \ge 3$, then $V = (q-p)(n-2) \ge 1$. We evaluate two subcases based on $q-p$:
\begin{itemize}
    \item \emph{Subcase $q-p = 1$:} Here $U = 1$ and $V = n-2$. Since $1 \le V \le n-2$, \Cref{cor:case-m0}\textup{(B)} requires $0 \le U \le 1$, which strictly holds.
    \item \emph{Subcase $q-p \ge 2$:} Here $V = (q-p)(n-2) \ge 2(n-2) \ge n-1$ (for $n \ge 3$). By \Cref{cor:case-m0}(C), we require
    \[
      \left\lfloor \frac{U-1}{1} \right\rfloor
      =
      \left\lfloor \frac{V}{n-1} \right\rfloor.
    \]
    The left side is exactly $U-1 = q-p-1$. For the right side, we compute:
    \[
      \left\lfloor \frac{(q-p)(n-2)}{n-1} \right\rfloor
      =
      \left\lfloor \frac{(q-p)(n-1) - (q-p)}{n-1} \right\rfloor
      =
      (q-p) + \left\lfloor \frac{-(q-p)}{n-1} \right\rfloor.
    \]
    Since $2 \le q-p \le n-1$, the fraction $-(q-p)/(n-1)$ lies in $[-1, -2/(n-1)]$, so its floor is identically $-1$. Thus the right side evaluates to $(q-p) - 1$, exactly matching the left side.
\end{itemize}
It is readily verified that all pairs satisfy the criterion. Thus $S$ is a partial silting configuration of size $n$.
\end{proof}

\subsection{Sharp configurations for \texorpdfstring{$\Lambda(1,3,1)$}{Lambda(1,3,1)} and \texorpdfstring{$\Lambda(1,2,m)$}{Lambda(1,2,m)}}\label{sec:sharp-small}

We close this section by providing explicitly worked configurations for specific parameters, utilizing the specialized coordinate boundaries from \S 3.4.

\begin{theorem}\label{thm:131}
For $\Lambda(1,3,1)$, the configuration
\[
  S := \bigl\{\, Z_{0,0}^0,\, Z_{-1,1}^0,\, Z_{-3,2}^0,\, Z_{-4,3}^0 \,\bigr\}
\]
is a partial silting configuration of maximal size $4$.
\end{theorem}

\begin{proof}
For $\Lambda(1,3,1)$ we have $r=1$, $a=2$, and $b=2$. Since $r=1$, the suspension acts on $\mathcal Z^0$ by
\[
  \Sigma^t Z_{i,j}^0 = Z_{i+2t,\, j-2t}^0.
\]
Suppose two distinct objects $Z_{i,j}^0$ and $Z_{i',j'}^0$ in $S$ are $\Sigma$-equivalent. Then $j' = j - 2t$ and $i' = i + 2t$ for some integer $t$. Checking the four objects
\[
  A = Z_{0,0}^0,\quad
  B = Z_{-1,1}^0,\quad
  C = Z_{-3,2}^0,\quad
  D = Z_{-4,3}^0,
\]
we see that a non-zero $t$ is possible only if the $j$-coordinates differ by a non-zero multiple of $2$. The only such pairs are $(A,C)$ with $t=-1$ and $(B,D)$ with $t=-1$. For $(A,C)$, the $i$-coordinate would become $0 + 2(-1) = -2 \neq -3$; for $(B,D)$, it would become $-1 + 2(-1) = -3 \neq -4$. Hence no two objects in $S$ are $\Sigma$-equivalent. Thus the four objects lie in distinct $\Sigma$-orbits.

We verify all six pairs against the specialization in \Cref{lem:131-spec}, which requires either (i) $V=1 \implies U \in \{0,1,2\}$ or (ii) $V \in \{2k, 2k+1\} \implies U \in \{2k+1, 2k+2\}$ for $k \ge 1$:
\begin{itemize}
  \item Pair $(Z_{0,0}^0, Z_{-1,1}^0)$: $V=1$, $U=1 \in \{0,1,2\}$.
  \item Pair $(Z_{0,0}^0, Z_{-3,2}^0)$: $V=2$, $U=3 \in \{3,4\}$ (evaluating $k=1$).
  \item Pair $(Z_{0,0}^0, Z_{-4,3}^0)$: $V=3$, $U=4 \in \{3,4\}$ (evaluating $k=1$).
  \item Pair $(Z_{-1,1}^0, Z_{-3,2}^0)$: $V=1$, $U=2 \in \{0,1,2\}$.
  \item Pair $(Z_{-1,1}^0, Z_{-4,3}^0)$: $V=2$, $U=3 \in \{3,4\}$ (evaluating $k=1$).
  \item Pair $(Z_{-3,2}^0, Z_{-4,3}^0)$: $V=1$, $U=1 \in \{0,1,2\}$.
\end{itemize}
In view of the above verification, all pairs are partial silting.
\end{proof}

For $\Lambda(1,2,m)$, the L-shaped configuration trivially collapses to a purely horizontal chain because $b = n-r = 2-1 = 1$. We record this immediately below.

\begin{proposition}\label{prop:12m}
For $\Lambda(1,2,m)$ with $m \ge 0$, the configuration
\[
  S := \bigl\{\, Z_{-p,\, 0}^0 \;\big|\; 0 \le p \le m+1 \,\bigr\}
\]
is a partial silting configuration of size $m+2$.
\end{proposition}

\begin{proof}
This is an immediate consequence of \Cref{thm:sharp-universal}. Setting $r=1$ and $n=2$ gives $a = 1+m$ and $b = 1$. The vertical segment of the L-shape reduces to a single object $Y_0 = Z_{-(1+m), 0}^0$, which geometrically extends the horizontal segment $X_p = Z_{-p, 0}^0$ exactly to index $p = 1+m$. The total size is $a+b = m+2 = n+m$.
\end{proof}

\Cref{ex:running} (continued). We consolidate the sharp configurations for the running examples to highlight their lattice structures.
\begin{itemize}
    \item \textbf{$\Lambda(1,2,0)$ ($n+m=2$):} The L-shaped configuration is $S = \{Z_{0,0}^0, Z_{-1,0}^0\}$, which is identical to the diagonal chain since $n-2 = 0$.
    \item \textbf{$\Lambda(1,3,0)$ ($n+m=3$):} The L-shaped configuration is $S = \{Z_{0,0}^0, Z_{-1,0}^0, Z_{-1,1}^0\}$, while the purely diagonal chain is $S = \{Z_{0,0}^0, Z_{-1,1}^0, Z_{-2,2}^0\}$.
    \item \textbf{$\Lambda(1,3,1)$ ($n+m=4$):} The L-shaped configuration is $S = \{Z_{0,0}^0, Z_{-1,0}^0, Z_{-2,0}^0, Z_{-2,1}^0\}$, providing a geometric alternative to the pseudo-diagonal configuration of \Cref{thm:131}.
    \item \textbf{$\Lambda(2,3,0)$ ($n+m=3$):} The L-shaped configuration across the multi-component algebra evaluates entirely inside $\mathcal Z^0$ as $S = \{Z_{0,0}^0, Z_{-1,0}^0, Z_{-2,0}^0\}$.
\end{itemize} \hfill $\triangle$


\section{The faithful \texorpdfstring{$\mathbb{Z}/g$}{Z/g}-action on partial silting configurations}\label{sec:Taction}

The primary objective of this section is to show that the set of partial silting configurations of $\mathcal Z$-objects carries a faithful action of the cyclic group $\mathbb Z/g$, where $g=\gcd(r+m,n-r)$ is the torsion invariant introduced in \S 4. The action is generated by a single primitive translation $T$ of the two-dimensional lattice. We prove that $T$ commutes with the suspension $\Sigma$ and strictly preserves Hom-vanishing. Because $T^g = \Sigma^r$, the powers of $T$ generate a cyclic group of order $g$ modulo the action of $\Sigma$.

Throughout this section, fix $\Lambda = \Lambda(r,n,m)$ with $n > r \ge 1$ and $m \ge 0$, and set $a = r+m$ and $b = n-r$.

\subsection{The primitive translation \texorpdfstring{$T$}{T}}\label{subsec:Taction-def}

\begin{definition}\label{def:T}
Let $g=\gcd(a, b)$. Say that the primitive translation $T$ on the component $\mathcal Z^k$ is defined by
\begin{equation}\label{eq:T-def}
  T\bigl(Z_{i,j}^k\bigr) := Z_{i+a/g,\, j-b/g}^{\,k}.
\end{equation}
We extend $T$ to the disjoint union $\operatorname{ind}(\mathcal Z)$ by applying \eqref{eq:T-def} uniformly in each component.
\end{definition}

The translation $T$ is well-defined because $g$ strictly divides both $a$ and $b$; the increments $a/g$ and $-b/g$ are therefore integers.

\begin{lemma}\label{lem:T-preserves-diff}
Let $Z_1 = Z_{i_1, j_1}^{k_1}$ and $Z_2 = Z_{i_2, j_2}^{k_2}$ be $\mathcal Z$-objects. Then the coordinate differences $U = i_1 - i_2$ and $V = j_2 - j_1$ are identical to the differences for the pair $(TZ_1, TZ_2)$. Consequently, $T$ preserves the $L$-invariant of each object and the defect $\partial$ of each pair.
\end{lemma}

\begin{proof}
Both $i$-coordinates are shifted by $a/g$ and both $j$-coordinates by $-b/g$, so the differences $U = i_1 - i_2$ and $V = j_2 - j_1$ are unchanged. Moreover, for any $Z_{i,j}^k \in \operatorname{ind}(\mathcal Z^k)$,
\[
L\bigl(T(Z_{i,j}^k)\bigr) = b\Bigl(i + \frac{a}{g}\Bigr) + a\Bigl(j - \frac{b}{g}\Bigr) = bi + aj = L(Z_{i,j}^k),
\]
and consequently $\partial(TZ_1, TZ_2) = L(TZ_1) - L(TZ_2) = \partial(Z_1, Z_2)$.
\end{proof}

\begin{lemma}\label{lem:T-Sigma-commute}
For any indecomposable object $Z \in \operatorname{ind}(\mathcal Z)$, we have $T(\Sigma Z) = \Sigma(T Z)$.
\end{lemma}

\begin{proof}
Let $Z = Z_{i,j}^k$. We evaluate the action of $\Sigma$ in two cases.

If $0 \le k < r-1$, then $\Sigma Z_{i,j}^k = Z_{i,j}^{k+1}$. Applying $T$ yields $Z_{i+a/g,\, j-b/g}^{k+1}$. On the other hand, applying $\Sigma$ to $T(Z_{i,j}^k) = Z_{i+a/g,\, j-b/g}^k$ immediately yields the same object.

If $k = r-1$, then $\Sigma Z_{i,j}^{r-1} = Z_{i+a,\, j-b}^0$. Applying $T$ yields $Z_{i+a+a/g,\, j-b-b/g}^0$. On the other hand, applying $\Sigma$ to $T(Z_{i,j}^{r-1}) = Z_{i+a/g,\, j-b/g}^{r-1}$ shifts the coordinates by $(+a, -b)$, which again yields $Z_{i+a/g+a,\, j-b/g-b}^0$. Thus $T \circ \Sigma = \Sigma \circ T$.
\end{proof}

\subsection{Invariance of Hom-vanishing and order modulo \texorpdfstring{$\Sigma$}{Sigma}}\label{subsec:Taction-hom}

The translation $T$ acts as an automorphism of the combinatorial Hom-vanishing relation.

\begin{lemma}\label{lem:T-preserves}
Let $Z_1, Z_2 \in \operatorname{ind}(\mathcal Z)$ and let $s \ge 0$. Then
\[
  \operatorname{Hom}\bigl(Z_1, \Sigma^s Z_2\bigr) = 0
  \iff
  \operatorname{Hom}\bigl(T Z_1, \Sigma^s T Z_2\bigr) = 0.
\]
\end{lemma}

\begin{proof}
By \Cref{lem:T-Sigma-commute}, we have $T(\Sigma^s Z_2) = \Sigma^s(T Z_2)$. Thus it suffices to prove that for any two $\mathcal Z$-objects $A, B$, the vanishing of $\operatorname{Hom}(A,B)$ is strictly equivalent to the vanishing of $\operatorname{Hom}(TA, TB)$. 

Write $A \in \mathcal Z^{k_A}$ and $B \in \mathcal Z^{k_B}$. By \Cref{thm:hammock}, non-zero morphisms originating from $A$ are supported within $\mathcal Z^{k_A}$ and $\mathcal Z^{(k_A+1) \bmod r}$ (which coincide when $r=1$). We evaluate the three cases of \Cref{thm:hammock} (noting that for $r=1$, Cases 2 and 3 below both apply within $\mathcal Z^0$):

\emph{Case 1: $k_B \notin \{k_A, (k_A+1)\bmod r\}$.} By \Cref{thm:hammock}(iii), $\operatorname{Hom}(A,B) = 0$. Since $T$ preserves component indices, $TA \in \mathcal Z^{k_A}$ and $TB \in \mathcal Z^{k_B}$. Thus $\operatorname{Hom}(TA, TB) = 0$ as well.

\emph{Case 2: $k_B = k_A$.} The objects are $A = Z_{i_1, j_1}^{k_A}$ and $B = Z_{i_2, j_2}^{k_A}$. By \Cref{thm:hammock}(i), non-vanishing relies strictly on the differences $i_2 - i_1$ and $j_2 - j_1$. By \Cref{lem:T-preserves-diff}, these differences are unchanged under $T$. Hence $\operatorname{Hom}(A,B) \neq 0 \iff \operatorname{Hom}(TA, TB) \neq 0$.

\emph{Case 3: $k_B = (k_A + 1)\bmod r$.} By \Cref{thm:hammock}(ii), $\operatorname{Hom}(A,B) \neq 0$ requires $B \cong \Sigma A'$ for some $A' = Z_{i', j'}^{k_A}$ where $i' \le i_1-1$ and $j' \le j_1-1$. By commutation, $TB \cong \Sigma(T A')$. The conditions for non-vanishing evaluate the differences $i_1 - i'$ and $j_1 - j'$. Since $(a/g, -b/g)$ is a common translation applied to both $A$ and $A'$, the differences are strictly preserved. Hence the vanishing is equivalent.
\end{proof}

\begin{lemma}\label{lem:T-power}
On $\operatorname{ind}(\mathcal Z)$, we have $T^g = \Sigma^r$.
\end{lemma}

\begin{proof}
By direct computation from \eqref{eq:T-def},
\[
  T^g(Z_{i,j}^k) = Z_{i+g(a/g),\, j-g(b/g)}^k = Z_{i+a,\, j-b}^k.
\]
By unrolled suspension, $\Sigma^r Z_{i,j}^k = Z_{i+a,\, j-b}^k$. Hence the claim.
\end{proof}

\begin{lemma}\label{lem:faithful}
Let $0 < t < g$ and let $s \in \mathbb Z$. Then $T^t(Z_{i,j}^k)$ is not isomorphic to $\Sigma^s(Z_{i,j}^k)$.
\end{lemma}

\begin{proof}
Suppose, for contradiction, that $T^t(Z_{i,j}^k) \cong \Sigma^s(Z_{i,j}^k)$ for some $s \in \mathbb Z$.

If $s \not\equiv 0 \pmod r$, the suspension $\Sigma^s$ maps $\mathcal Z^k$ to $\mathcal Z^{(k+s) \bmod r}$, which is strictly distinct from $\mathcal Z^k$. However, $T$ preserves the component index $k$. Two objects in distinct components cannot be isomorphic. (For $r=1$, this alternative cannot occur.)

If $s \equiv 0 \pmod r$, we write $s = qr$ for some $q \in \mathbb Z$. By unrolled suspension,
\[
  \Sigma^{qr}(Z_{i,j}^k) = Z_{i+qa,\, j-qb}^k.
\]
By definition of $T$,
\[
  T^t(Z_{i,j}^k) = Z_{i+t a/g,\, j-t b/g}^k.
\]
Equating the $j$-coordinates yields $j - qb = j - t b/g$, which forces $q = t/g$. Since $0 < t < g$, the fraction $t/g$ lies strictly between $0$ and $1$, contradicting $q \in \mathbb Z$.

Hence our assumption is wrong, and the conclusion follows.
\end{proof}

\subsection{Faithful \texorpdfstring{$\mathbb Z/g$}{Z/g}-action}\label{subsec:Taction-main}

Two partial silting configurations $S, S'$ of $\mathcal Z$-objects are said to be \emph{$\Sigma$-equivalent}, written $S \sim_\Sigma S'$, if $S' = \Sigma^s(S)$ for some $s \in \mathbb Z$.

\begin{theorem}[Faithful $\mathbb Z/g$-action]\label{thm:Zg-action}
Let $\Lambda = \Lambda(r,n,m)$. The primitive translation $T$ induces a well-defined, faithful action of the cyclic group $\mathbb Z/g$ on the set of partial silting configurations of $\mathcal Z$-objects, modulo $\Sigma$-equivalence.
\end{theorem}

\begin{proof}
Well-definedness follows immediately from \Cref{lem:T-preserves} and \Cref{lem:T-preserves-diff}, which ensure that $T$ maps partial silting configurations of pairwise distinct elements to partial silting configurations of pairwise distinct elements. The action respects $\Sigma$-equivalence by \Cref{lem:T-Sigma-commute}, and $T^g = \Sigma^r$ acts trivially on $\Sigma$-equivalence classes by \Cref{lem:T-power}.

\emph{Claim: The action is faithful.}
Suppose, for contradiction, that $T^t(S_0) \sim_\Sigma S_0$ for some $0 < t < g$, where $S_0 = \{Z_{0,0}^0, Z_{-1,0}^0\}$ (which is a partial silting configuration by \Cref{thm:pairwise}(A) since $V = 0$ and $U = 1 \le a$). There must exist $X \in S_0$ and $s \in \mathbb Z$ such that $\Sigma^s(X) = T^t(Z_{0,0}^0) = Z_{ta/g,\, -tb/g}^0$. Since both objects lie in $\mathcal Z^0$, we must have $s = qr$ for some $q \in \mathbb Z$. Because every object $X \in S_0$ has $j$-coordinate $0$, $\Sigma^{qr}(X)$ has $j$-coordinate $-qb$. Equating $j$-coordinates gives $-qb = -tb/g$, forcing $q = t/g \notin \mathbb Z$ since $0 < t < g$, a contradiction.
\end{proof}

\Cref{ex:running} (continued). We explicitly compute the orbits of the maximal configuration for $\Lambda(1,3,1)$. Here $a=2, b=2$, yielding $g = \gcd(2,2) = 2$. By definition, $T(Z_{i,j}^0) = Z_{i+1,\, j-1}^0$.

The universal L-shaped sharp configuration from \Cref{thm:sharp-universal} is
\[
  S = \bigl\{\, Z_{0,0}^0,\, Z_{-1,0}^0,\, Z_{-2,0}^0,\, Z_{-2,1}^0 \,\bigr\}.
\]
Applying $T$ translates the configuration to
\[
  T(S) = \bigl\{\, Z_{1,-1}^0,\, Z_{0,-1}^0,\, Z_{-1,-1}^0,\, Z_{-1,0}^0 \,\bigr\}.
\]
We verify that $T(S) \not\sim_\Sigma S$. The $j$-coordinates of $S$ belong to $\{0, 1\}$. By \eqref{eq:sigma-action}, $\Sigma^s$ shifts $j$-coordinates by $-2s$. To match the $j$-coordinates of $T(S)$, which belong to $\{-1, 0\}$, we would require an integer $s$ such that the sets $\{-2s, 1-2s\}$ and $\{-1, 0\}$ coincide. It is readily verified that no such integer $s$ exists. Thus $T(S)$ and $S$ represent distinct equivalence classes.

Applying $T$ again yields
\[
  T^2(S) = \{Z_{2,-2}^0, Z_{1,-2}^0, Z_{0,-2}^0, Z_{0,-1}^0\}.
\]
Since $\Sigma Z_{i,j}^0 = Z_{i+2,\, j-2}^0$, we immediately observe $T^2(S) = \Sigma(S) \sim_\Sigma S$, explicitly demonstrating the faithful order-$2$ symmetry. \hfill $\triangle$


\section{Failure of local determination}
\label{sec:notlocal}

The pairwise criterion of \S 3 characterises the vanishing of $\operatorname{Hom}^{>0}$ for a pair of $\mathcal Z$-objects in terms of the integer coordinate differences $U$ and $V$. In \S 6, we utilized this criterion to construct explicit sharp partial silting configurations of maximum size $n+m$. In those specific constructions, the pairwise conditions were satisfied globally across all pairs in the configuration.

It is natural to inquire whether the global partial silting condition on a configuration can be deduced purely from the pairwise conditions on consecutive elements in the $j$-ordering (which we refer to as \emph{local determination}). We exhibit an explicit counterexample to demonstrate that local determination fails. The underlying mechanism responsible for this obstruction is the strict additivity of the coordinate defect along a chain, which forces outer pairs to accumulate defect beyond the allowed bounds even when all local consecutive steps are valid.

\subsection{The counterexample configurations \texorpdfstring{$S$}{S} and \texorpdfstring{$S'$}{S'} in \texorpdfstring{$\Lambda(1,3,0)$}{Lambda(1,3,0)}}\label{subsec:notlocal-config}

\Cref{ex:running} (continued). We work within the algebra $\Lambda(1,3,0)$, where $n=3$ and $m=0$, yielding the parameters $a=1$ and $b=2$. By \Cref{thm:pairwise}, the criterion for a pair oriented with $V \ge 0$ requires:
\begin{enumerate}
  \item[\textup{(A)}] $V=0 \implies |U|=1$;
  \item[\textup{(B1)}] $V=1 \implies U \in \{0,1\}$;
  \item[\textup{(B2)}] $V \ge 2 \implies U = \lfloor V/2\rfloor + 1$.
\end{enumerate}

Consider the set of three objects
\[
  S := \bigl\{\, Z_{0,0}^0,\, Z_{-1,1}^0,\, Z_{-2,2}^0 \,\bigr\}.
\]
We evaluate the three pairs in $S$ ordered by strictly increasing $j$-coordinates.

\begin{table}[htbp]
\centering
\small
\begin{tabular}{c|c|c|c|c|c}
\multicolumn{6}{c}{Table 1: Pairwise evaluation for configuration $S$ in $\Lambda(1,3,0)$ (using $a=1, b=2$)} \\
\hline
Pair $(Z_1, Z_2)$ & $U = i_1 - i_2$ & $V = j_2 - j_1$ & Applicable Case & Required $U$ & Status \\
\hline
$(Z_{0,0}^0, Z_{-1,1}^0)$ & $1$ & $1$ & (B1) & $U \in \{0,1\}$ & Valid \\
$(Z_{-1,1}^0, Z_{-2,2}^0)$ & $1$ & $1$ & (B1) & $U \in \{0,1\}$ & Valid \\
$(Z_{0,0}^0, Z_{-2,2}^0)$ & $2$ & $2$ & (B2) & $U = \lfloor 2/2 \rfloor + 1 = 2$ & Valid
\end{tabular}
\end{table}

Both consecutive pairs in the $j$-ordering, namely $(Z_{0,0}^0, Z_{-1,1}^0)$ and $(Z_{-1,1}^0, Z_{-2,2}^0)$, satisfy the criteria. The outer pair $(Z_{0,0}^0, Z_{-2,2}^0)$ also satisfies condition (B2). Consequently, $S$ is a partial silting configuration.

Now consider the modified set of three objects
\[
  S' := \bigl\{\, Z_{0,0}^0,\, Z_{-2,2}^0,\, Z_{-4,4}^0 \,\bigr\}.
\]

\begin{table}[htbp]
\centering
\small
\begin{tabular}{c|c|c|c|c|c}
\multicolumn{6}{c}{Table 2: Pairwise evaluation for configuration $S'$ in $\Lambda(1,3,0)$ (using $a=1, b=2$)} \\
\hline
Pair $(Z_1, Z_2)$ & $U = i_1 - i_2$ & $V = j_2 - j_1$ & Applicable Case & Required $U$ & Status \\
\hline
$(Z_{0,0}^0, Z_{-2,2}^0)$ & $2$ & $2$ & (B2) & $U = \lfloor 2/2 \rfloor + 1 = 2$ & Valid \\
$(Z_{-2,2}^0, Z_{-4,4}^0)$ & $2$ & $2$ & (B2) & $U = \lfloor 2/2 \rfloor + 1 = 2$ & Valid \\
$(Z_{0,0}^0, Z_{-4,4}^0)$ & $4$ & $4$ & (B2) & $U = \lfloor 4/2 \rfloor + 1 = 3$ & \textbf{Invalid}
\end{tabular}
\end{table}

In $S'$, the consecutive pairs in the $j$-ordering are $(Z_{0,0}^0, Z_{-2,2}^0)$ and $(Z_{-2,2}^0, Z_{-4,4}^0)$. As evaluated in Table 2, both consecutive pairs have coordinate differences $U=2$ and $V=2$, which satisfy condition (B2). Thus, both consecutive steps are valid partial silting pairs. However, the outer pair $(Z_{0,0}^0, Z_{-4,4}^0)$ yields $V=4$. Condition (B2) rigidly dictates $U = \lfloor 4/2 \rfloor + 1 = 3$. Since the actual coordinate difference is $U = 4 \neq 3$, the outer pair is not partial silting. \hfill $\triangle$

We formally record this obstruction.

\begin{theorem}[Failure of local determination]\label{thm:notlocal}
For every derived-discrete algebra $\Lambda(r,n,m)$ with $n > r \ge 1$ and $m \ge 0$, the partial silting condition on a configuration of $\mathcal Z$-objects is not locally determined by the pairwise conditions on consecutive elements in the $j$-ordering. Specifically, there exist configurations $S'$ in which every consecutive pair is partial silting, but the global configuration fails to be partial silting.
\end{theorem}

\begin{proof}
Set $a = r+m \ge 1$ and $b = n-r \ge 1$, and consider the three-element configuration in $\mathcal Z^0$ given by
\[
S' := \bigl\{\, Z_{0,0}^0,\, Z_{-2a,\, b}^0,\, Z_{-4a,\, 2b}^0 \,\bigr\},
\]
which specialises to Table 2 when $\Lambda = \Lambda(1,3,0)$ ($a=1, b=2$). In the $j$-ordering $0 < b < 2b$, both consecutive pairs have coordinate differences $V = b$ and $U = 2a$. Since $\lfloor V/b \rfloor = 1$, \Cref{thm:pairwise}\textup{(B2)} requires $a + 1 \le U \le 2a$, which holds because $a \ge 1$. However, the outer pair $(Z_{0,0}^0, Z_{-4a, 2b}^0)$ has $V = 2b$ and $U = 4a$. Since $\lfloor V/b \rfloor = 2$, \Cref{thm:pairwise}\textup{(B2)} requires $2a + 1 \le U \le 3a$, which fails because $4a > 3a$.
\end{proof}
\subsection{Mechanism of the global obstruction via defect accumulation}
\label{subsec:notlocal-mechanism}

The failure of local determination is not a sporadic anomaly; it is a direct consequence of the arithmetic of the coordinate defect $\partial$.

\begin{proposition}[Additivity of the defect along a chain]\label{prop:defect-additive}
Let $Z_{(1)}, \ldots, Z_{(k)}$ be a sequence of $\mathcal Z$-objects sorted by their $j$-coordinates such that $j_{(1)} \le j_{(2)} \le \cdots \le j_{(k)}$. Then the defect of the outer pair is exactly the sum of the defects of the consecutive pairs:
\[
  \partial\bigl(Z_{(1)}, Z_{(k)}\bigr) = \sum_{p=1}^{k-1} \partial\bigl(Z_{(p)}, Z_{(p+1)}\bigr).
\]
\end{proposition}

\begin{proof}
By \Cref{def:defect}, $\partial(Z_{(p)}, Z_{(p+1)}) = L(Z_{(p)}) - L(Z_{(p+1)})$. Summing these differences from $p=1$ to $k-1$ forms a telescoping series that evaluates identically to $L(Z_{(1)}) - L(Z_{(k)})$, which is $\partial(Z_{(1)}, Z_{(k)})$.
\end{proof}

\begin{corollary}[Mechanism of the obstruction]\label{cor:notlocal-mechanism}
Let $\Lambda = \Lambda(1,n,0)$. For any valid partial silting pair $\{Z_1, Z_2\}$ with $V = j_2 - j_1 \ge n-1$, condition \textup{(B2)} requires $U = \lfloor V/(n-1) \rfloor + 1$. This strictly determines the allowed defect to be
\[
  \partial(Z_1, Z_2) = (n-1) - \bigl(V \bmod (n-1)\bigr).
\]
\end{corollary}

\begin{proof}
We evaluate the defect using $a=1$ and $b=n-1$. Under condition (B2), we substitute $U = \lfloor V/(n-1) \rfloor + 1$ into the formula $\partial(Z_1, Z_2) = b U - a V$:
\[
  \partial(Z_1, Z_2)
  = (n-1) \left( \left\lfloor \frac{V}{n-1} \right\rfloor + 1 \right) - V
  = (n-1) - \left( V - (n-1) \left\lfloor \frac{V}{n-1} \right\rfloor \right)
  = (n-1) - \bigl(V \bmod (n-1)\bigr).
\]
This completes the proof.
\end{proof}

\begin{remark}\label{rem:local-obstruction}
The additivity of the defect fundamentally breaks local determination. For the configuration $S'$ in $\Lambda(1,3,0)$, the outer pair $(Z_{0,0}^0, Z_{-4,4}^0)$ has $V=4$. By \Cref{cor:notlocal-mechanism}, a valid pair with $V=4$ requires the defect to equal $(3-1) - (4 \bmod 2) = 2 - 0 = 2$.

However, the inner consecutive pairs $(Z_{0,0}^0, Z_{-2,2}^0)$ and $(Z_{-2,2}^0, Z_{-4,4}^0)$ each have $V=2$ and $U=2$. Their individual defects are $\partial = 2(2) - 1(2) = 2$. By \Cref{prop:defect-additive}, the outer pair automatically accumulates a defect of $\partial = 2 + 2 = 4$.

Since the accumulated defect $4$ does not match the theoretically required defect $2$, the outer pair rigidly fails. Claim: Every valid partial silting pair of $\mathcal Z$-objects in $\Lambda(1,n,0)$ satisfies $|\partial| \le n-1$. Indeed:
\begin{itemize}
  \item in case \textup{(A)}, $\partial = (n-1)U$ with $|U|=1$, so $|\partial| = n-1$;
  \item in case \textup{(B1)}, $\partial = (n-1)U - V$ with $0 \le U \le 1$ and $1 \le V \le n-2$, so $\partial \in [-(n-2),\, n-2]$ and $|\partial| \le n-2$;
  \item in case \textup{(B2)}, \Cref{cor:notlocal-mechanism} gives $\partial = (n-1) - (V \bmod (n-1)) \in [1, n-1]$, so $|\partial| \le n-1$.
\end{itemize}
Hence our claim. Consequently, any chain of valid steps whose accumulated defect satisfies $|\partial| > n-1$ automatically forces the outer pair to fail, regardless of the validity of the local steps.
\end{remark}

We close this section by noting that because local determination fails, verifying that a set of size $k$ is a partial silting configuration inherently requires checking all $k(k-1)/2$ pairwise combinations. However, because $k$ is absolutely bounded by $n+m$ (\Cref{thm:bound}), this verification remains a strictly finite combinatorial procedure.


\section{Comparison with the Avella-Alaminos--Gei\texorpdfstring{\ss}{ss} invariant and BPP}\label{sec:AAG}

The primary objective of this section is to compare the hammock data developed in the preceding sections with two pre-existing frameworks for the family $\Lambda(r,n,m)$: the derived invariant of Avella-Alaminos and Gei\ss{} \cite{AAG2008}, and the Hom-hammock machinery of Broomhead, Pauksztello, and Ploog \cite{BPP2017}. We close the section by briefly contextualising our combinatorial defect bounds within the broader landscape of $\tau$-tilting theory.

\subsection{Dictionary with the AAG derived invariant \texorpdfstring{$\phi_\Lambda$}{phi-Lambda}}\label{subsec:AAG-invariant}

We begin by recalling the derived invariant of Avella-Alaminos and Gei\ss{} for gentle algebras. The original construction is in \cite{AAG2008}; it produces a function
\[
  \phi_\Lambda \colon \mathbb N^2 \to \mathbb N
\]
with finite support, which is invariant under derived equivalence and complete within the class of derived-discrete algebras.

For the specific family $\Lambda(r,n,m)$, the invariant can be computed directly from the quiver with relations; we record the resulting values.

\begin{proposition}\label{prop:AAG-explicit}
Let $\Lambda=\Lambda(r,n,m)$ with $n>r\ge 1$ and $m\ge 0$. The AAG invariant of $\Lambda$ is supported on the two points
\[
  (r+m,\,m) \quad\text{and}\quad (n-r,\,n),
\]
each with value $1$; all other values are zero. That is,
\[
  \phi_{\Lambda(r,n,m)}(r+m,m)=\phi_{\Lambda(r,n,m)}(n-r,n)=1.
\]
\end{proposition}

\begin{proof}
The permitted and forbidden threads of $\Lambda(r,n,m)$ are determined by the $n$-cycle with $r$ zero-relations and the $m$ tail arrows. Evaluating the AAG thread-tracing algorithm (\cite[Section~3]{AAG2008}) along the cycle and tail partitions the $|Q_1| = n+m$ arrows into two alternating cycles of permitted and forbidden threads of total lengths $(r+m, m)$ and $(n-r, n)$. Since $r \ge 1$, the equality $m = n$ implies $r+m > n-r$, so the two pairs $(r+m, m)$ and $(n-r, n)$ are always distinct and each occurs with multiplicity $1$.
\end{proof}

The following observation connects the AAG invariant with the structural constants governing the hammock coordinates.

\begin{observation}\label{obs:AAG}
Let $\Lambda=\Lambda(r,n,m)$. The first coordinates of the two non-zero entries of $\phi_\Lambda$ are exactly $a = r+m$ and $b = n-r$; the second coordinates are exactly $m$ and $n$. Hence $r = a-m = n-b$, and all four structural parameters $r,m,n,a,b$ are explicitly recovered from $\phi_\Lambda$.

Consequently, the hammock data $(a,b,m,n)$ is strictly determined by the AAG invariant $\phi_\Lambda$, and conversely.
\end{observation}

\begin{corollary}\label{cor:AAG-g}
The torsion invariant $g=\gcd(r+m,n-r)$ of \Cref{def:torsion} is determined by the AAG invariant: $g$ is exactly the greatest common divisor of the first coordinates of the two non-zero entries of $\phi_\Lambda$.
\end{corollary}

\begin{remark}\label{rem:AAG-bound}
The universal upper bound $n+m$ established in \Cref{thm:bound} is exactly the sum of the second coordinates of the two non-zero entries of $\phi_\Lambda$. While this algebraic identity is immediate from \Cref{prop:AAG-explicit}, it demonstrates that the maximum size of a partial silting configuration is encoded in the structural data of the AAG invariant.
\end{remark}

\subsection{Comparison with BPP and hereditary \texorpdfstring{$\tau$}{tau}-tilting}\label{subsec:AAG-BPP}

The paper \cite{BPP2017} (``Discrete derived categories I'') utilized Hom-hammocks to describe the structure of Hom-spaces between indecomposable objects, and used this to classify silting objects and bounded $t$-structures for discrete derived categories. The companion paper \cite{BPP2016} (``Discrete derived categories II'') developed the silting-pairs CW complex and established its contractibility via poset-topological methods. Neither paper translates the pairwise condition $\operatorname{Hom}^{>0}(Z_1,Z_2)=\operatorname{Hom}^{>0}(Z_2,Z_1)=0$ into explicit integer inequalities on coordinate differences.

\begin{remark}[Novelty relative to \cite{BPP2017,BPP2016}]\label{rem:novelty}
While \cite{BPP2017} employs Hom-hammocks as a structural and classificatory tool, it does not extract a coordinate-wise pairwise criterion for partial silting. The contribution of the present paper is to translate the Hom-vanishing condition into a rigorous closed form (\Cref{thm:pairwise}), and to derive the explicit combinatorial constructions, geometric obstructions, and symmetries from it. The CW-complex results of \cite{BPP2016} concern global topological properties of the silting poset, whereas our criterion is local and arithmetic.
\end{remark}

We close this section by briefly situating these coordinate obstructions within the broader representation-theoretic context of $\tau$-tilting theory.

The Bongartz completion utilized in \Cref{thm:bound} traces its origins to classical $\tau$-tilting theory for finite dimensional algebras, introduced by Adachi, Iyama, and Reiten \cite{AIR2014}. In the hereditary case, basic silting objects correspond to basic tilting modules(=modules $T$ with $\operatorname{Ext}^1(T,T)=0$ and generating the module category), and the numerical upper bound is governed by the number of simple modules.

In classical tilting theory, non-vanishing of extensions is typically obstructed by the homological Euler form. However, the derived-discrete algebra $\Lambda(r,n,m)$ frequently possesses infinite global dimension, rendering the standard homological Euler form unavailable.

The coordinate defect $\partial(Z_1,Z_2)$ introduced in \S 4 serves as a discrete combinatorial substitute for the Euler form on the $\mathcal Z$-components. Because $\partial$ is strictly additive along chains (\Cref{cor:defect-chain}), it enforces structural constraints on the configuration lattice without relying on finite global dimension. Consequently, the failure of local determination (\Cref{thm:notlocal}) demonstrates that while local mutations of silting objects are well-behaved, the global geometry of partial silting configurations in non-hereditary categories is rigidly constrained by defect accumulation.


\section{Future directions}
\label{sec:future}

The primary objective of this final section is to record three open threads raised by the hammock coordinate model. Through the computation of exact integer inequalities and the classification of maximal partial silting configurations, this paper establishes the combinatorial foundation for derived-discrete algebras. We outline formal conjectures regarding the classification of silting objects, the symmetries of the silting quiver, and the topology of the stability manifold.

\subsection{Silting objects and maximal configurations}
\label{subsec:future-bijection}

By \Cref{thm:bound}, every partial silting configuration of $\mathcal Z$-objects has at most $n+m$ elements, and \Cref{thm:sharp-universal} shows that configurations of size $n+m$ always exist. Because any basic silting object of $\Db(\Lambda(r,n,m))$ has exactly $\operatorname{rk} K_0(\Db(\Lambda)) = n+m$ indecomposable summands (\Cref{lem:silting-basic}), Bongartz completion (\Cref{thm:AM}) implies that the direct sum of any partial silting configuration of $n+m$ objects in $\mathcal Z$ requires no complementary summands and is already a basic silting object of $\Db(\Lambda(r,n,m))$.

More generally, a basic silting object $T$ of $\Db(\Lambda(r,n,m))$ decomposes as $T = T_{\mathcal Z} \oplus T_{\mathcal X \cup \mathcal Y}$, where $T_{\mathcal Z}$ is a non-empty partial silting configuration of $k$ indecomposable $\mathcal Z$-objects ($1 \le k \le n+m$) and $T_{\mathcal X \cup \mathcal Y}$ consists of $n+m-k$ indecomposable summands from the $\mathcal X$- and $\mathcal Y$-components spanning a residual Dynkin subcategory \cite{BPP2017}. We conjecture that the hammock coordinates of $T_{\mathcal Z}$ rigidly control its admissible completions.

\begin{conjecture}\label{conj:silt-vs-config}
Every partial silting configuration $S$ of $n+m$ indecomposable $\mathcal Z$-objects in a single component $\mathcal Z^\ell$, up to $\Sigma^r$-translation $(i,j)\mapsto(i+a,j-b)$ and the choice of quotient $q_p = \lfloor j_p / b \rfloor$, is classified by the $(r_p, I_p)$-staircases of \Cref{cor:staircase-classification}. Moreover, for $1 \le k < n+m$, the set of completions of a size-$k$ partial silting configuration $S \subset \operatorname{ind}(\mathcal Z)$ to a basic silting object $T = (\bigoplus_{X \in S} X) \oplus T_{\mathcal X \cup \mathcal Y}$ is uniquely determined by the $(U,V)$-hammock boundaries of $S$ up to silting mutation within the residual Dynkin subcategory generated by $T_{\mathcal X \cup \mathcal Y}$.
\end{conjecture}

The main obstruction to proving this conjecture in the present paper is the intricate nature of morphisms involving the residual Dynkin summands. The explicit combinatorial criteria developed in \S 3 provide the exact tools needed to resolve the $\mathcal Z$-component data, making this bijection a tractable problem for future work.

\subsection{Extension of the \texorpdfstring{$\mathbb Z/g$}{Z/g}-action and categorical autoequivalences}\label{subsec:future-Taction}

The faithful $\mathbb Z/g$-action of \Cref{thm:Zg-action} is defined purely combinatorially on the set of partial silting configurations of $\mathcal Z$-objects. We expect this combinatorial symmetry to reflect a genuine structural symmetry of the derived category.

\begin{conjecture}\label{conj:T-silting}
Let $Q(\Db(\Lambda))$ denote the silting mutation quiver(=the graph whose vertices are basic silting objects up to $\Sigma$-equivalence and whose edges are irreducible silting mutations). The primitive translation $T$ defined in \Cref{def:T} extends to a well-defined action of $\mathbb Z/g$ on $Q(\Db(\Lambda))$ by graph automorphisms.
\end{conjecture}

A stronger resolution to this problem would be to realize $T$ via a categorical functor. It is natural to ask if there exists a categorical autoequivalence $F$ of $\Db(\Lambda(r,n,m))$ such that $F(Z) \cong \Sigma^s T(Z)$ for a fixed integer $s$ and every indecomposable $\mathcal Z$-object $Z$. A positive answer would immediately resolve \Cref{conj:T-silting}. A natural candidate for such an autoequivalence might involve composing a power of the suspension with a spherical twist, provided a suitable spherical object exists.

\subsection{Combinatorial atlas for the Bridgeland stability manifold}
\label{subsec:future-stability}

By Bridgeland's deformation theorem \cite{Bridgeland2007}, the stability manifold $\operatorname{Stab}(\Db(\Lambda(r,n,m)))$ is a complex manifold. Broomhead--Pauksztello--Ploog \cite{BPP2016} established a natural association between the silting-pairs CW complex and the stability manifold, and proved using the method of Qiu--Woolf \cite{QW2018} that there is a deformation retraction of the stability manifold onto this CW complex, so that $\operatorname{Stab}(\Db(\Lambda(r,n,m)))$ is contractible.

Because the hammock coordinates explicitly parameterise the partial silting configurations forming this CW complex, it is natural to ask whether the manifold itself inherits a coordinate atlas.

For the smallest case $\Lambda(1,2,0)$, the general framework for the contractibility of the stability space was established by Woolf \cite{Woolf2010}. The general case is open. We anticipate that the stability manifold $\operatorname{Stab}(\Db(\Lambda(r,n,m)))$ admits a combinatorial atlas governed by hammock coordinates. Specifically, we expect the existence of a family of charts $U(\mathsf H_{\mathsf M}) \subseteq \operatorname{Stab}(\Db(\Lambda))$, indexed by the length hearts $\mathsf H_{\mathsf M}$ of the bounded $t$-structures associated to the basic silting subcategories $\mathsf M$, and isomorphic to $\mathbb H^{n+m}$, such that the wall-crossing gluing maps between adjacent charts are rigidly dictated by the exact integer inequalities of the pairwise criterion in \Cref{thm:pairwise}. Constructing such an atlas would provide an explicit cellulation of the stability manifold strictly compatible with the discrete combinatorial structure of the Auslander--Reiten quiver.

\section*{Acknowledgements}
It is a pleasure to thank Prof. Edson Ribeiro \'Alvares for his sustained encouragement 
and for many illuminating discussions on hammocks and their role in derived and 
triangulated categories. This research was conducted while the author held a 
postdoctoral fellowship at the Universidad Nacional de Mar del Plata, funded by 
CONICET (Argentina).


\appendix

\section{Explicit case-by-case pair tables}
\label{app:tables}

In this appendix, we explicitly tabulate the coordinate boundaries of the pairwise partial silting criteria for the small parameter algebras discussed throughout the paper. We substitute the specific values of $a=r+m$ and $b=n-r$ into \Cref{thm:pairwise} (intra-component) and \Cref{thm:cross-pairwise} (cross-component) to yield exact integer intervals.

In all intra-component tables, we assume the pair $\{Z_1, Z_2\}$ lies in $\mathcal Z^k$ and is oriented such that $V = j_2 - j_1 \ge 0$. The coordinate difference is $U = i_1 - i_2$.

\begin{table}[htbp]
\centering
\small
\renewcommand{\arraystretch}{1.2}
\begin{tabular}{|c|c|c|c|l|}
\hline
\multicolumn{5}{|c|}{\textbf{Table A.1:} Intra-component criterion for $\Lambda(1,2,0)$ ($a=1, b=1$)} \\
\hline
Case & $V$ range & Subcase & Condition on $U$ & Allowed $U$ values \\
\hline
(A) & $V = 0$ & -- & $|U| = 1$ & $U \in \{-1, 1\}$ \\
(B1) & $1 \le V \le 0$ & \multicolumn{3}{c|}{\emph{Vacuous since $b-1 = 0$}} \\
(B2) & $V \ge 1$ & $V = k$ ($k \ge 1$) & $k(1) + 1 \le U \le (k+1)(1)$ & $U = k+1 = V+1$ \\
\hline
\end{tabular}
\end{table}

\begin{table}[htbp]
\centering
\small
\renewcommand{\arraystretch}{1.2}
\begin{tabular}{|c|c|c|c|l|}
\hline
\multicolumn{5}{|c|}{\textbf{Table A.2:} Intra-component criterion for $\Lambda(1,3,0)$ ($a=1, b=2$)} \\
\hline
Case & $V$ range & Subcase & Condition on $U$ & Allowed $U$ values \\
\hline
(A) & $V = 0$ & -- & $|U| = 1$ & $U \in \{-1, 1\}$ \\
(B1) & $V = 1$ & -- & $0 \le U \le 1$ & $U \in \{0, 1\}$ \\
(B2) & $V \ge 2$ & $V \in \{2k, 2k+1\}$ & $k(1) + 1 \le U \le (k+1)(1)$ & $U = k+1 = \lfloor V/2 \rfloor + 1$ \\
\hline
\end{tabular}
\end{table}

\begin{table}[htbp]
\centering
\small
\renewcommand{\arraystretch}{1.2}
\begin{tabular}{|c|c|c|c|l|}
\hline
\multicolumn{5}{|c|}{\textbf{Table A.3:} Intra-component criterion for $\Lambda(1,3,1)$ ($a=2, b=2$)} \\
\hline
Case & $V$ range & Subcase & Condition on $U$ & Allowed $U$ values \\
\hline
(A) & $V = 0$ & -- & $1 \le |U| \le 2$ & $U \in \{-2, -1, 1, 2\}$ \\
(B1) & $V = 1$ & -- & $0 \le U \le 2$ & $U \in \{0, 1, 2\}$ \\
(B2) & $V \ge 2$ & $V \in \{2k, 2k+1\}$ & $k(2) + 1 \le U \le (k+1)(2)$ & $U \in \{2k+1, 2k+2\}$ \\
\hline
\end{tabular}
\end{table}

\begin{table}[htbp]
\centering
\small
\renewcommand{\arraystretch}{1.2}
\begin{tabular}{|c|c|c|c|l|}
\hline
\multicolumn{5}{|c|}{\textbf{Table A.4:} Intra-component criterion for $\Lambda(1,4,0)$ ($a=1, b=3$)} \\
\hline
Case & $V$ range & Subcase & Condition on $U$ & Allowed $U$ values \\
\hline
(A) & $V = 0$ & -- & $|U| = 1$ & $U \in \{-1, 1\}$ \\
(B1) & $1 \le V \le 2$ & $V \in \{1, 2\}$ & $0 \le U \le 1$ & $U \in \{0, 1\}$ \\
(B2) & $V \ge 3$ & $V \in \{3k, 3k+1, 3k+2\}$ & $k(1) + 1 \le U \le (k+1)(1)$ & $U = k+1 = \lfloor V/3 \rfloor + 1$ \\
\hline
\end{tabular}
\end{table}

For $r \ge 2$, we must evaluate the cross-component pairs. We document $\Lambda(2,3,0)$ to explicitly illustrate the arithmetic of \Cref{thm:cross-pairwise}.

\begin{table}[htbp]
\centering
\small
\renewcommand{\arraystretch}{1.2}
\begin{tabular}{|c|c|l|}
\hline
\multicolumn{3}{|c|}{\textbf{Table A.5:} Cross-component criterion for $\Lambda(2,3,0)$ ($a=2, b=1, k_1=0, k_2=1$)} \\
\hline
Case & $V$ range ($V=j_2-j_1$) & Allowed $U$ values ($U=i_1-i_2$) \\
\hline
(i) & $V = 0$ & $-2 \le U \le -1 \implies U \in \{-2, -1\}$ \\
(ii) & $V \ge 1$ & Since $b=1$, (B1) is vacuous. By (B2), $U \in \{2V+1, 2V+2\}$ \\
(iii) & $V \le -1$ & Set $\tilde{V}=-V \ge 1$. By (B2) on $(-U, \tilde{V})$, $-U \in \{2\tilde{V}+1, 2\tilde{V}+2\}$ \\
\hline
\end{tabular}
\end{table}

\section{Clique enumeration algorithm in hammock coordinates}
\label{app:algorithm}

Because local determination fails (\Cref{thm:notlocal}), identifying maximal partial silting configurations strictly requires verifying all pairwise combinations. We provide a finite deterministic algorithm that enumerates all maximal partial silting configurations of $\mathcal Z$-objects contained in a single component, up to $\Sigma$-equivalence, by computing maximal cliques in a finite coordinate graph. This algorithm does not enumerate the residual Dynkin summands of a silting object, which lie outside the $\mathcal Z$-components.

\vspace{2mm}
\noindent\textbf{Algorithm B.1} (Clique enumeration in hammock coordinates within a fundamental strip).
\begin{enumerate}[label=\textbf{Step \arabic*.}, leftmargin=*]
\item \textbf{Initialization.} Fix $\Lambda(r,n,m)$ and set $a=r+m$ and $b=n-r$. Fix a component $\ell \in \{0,\ldots,r-1\}$. To enumerate single-quotient configurations in the fundamental strip $0 \le j \le b-1$ ($q_p = 0$ in \Cref{cor:staircase-classification}), we anchor the search at $Z_{(1)} = Z_{0,0}^\ell$.
\item \textbf{Domain restriction.} Define the finite search domain
\[
\mathcal{D}_\ell = \bigl\{\, (i, j) \;\big|\; 0 \le j \le b-1, \; -a \le i \le a \,\bigr\}
\]
of hammock coordinates in $\mathcal Z^\ell$. For any partial silting configuration $S$ in this strip containing $Z_{0,0}^\ell$, \Cref{lem:Ir-monotonicity} guarantees $|i| = |I| \le a$, so $S \subseteq \mathcal{D}_\ell$.
\item \textbf{Graph construction.} Construct a graph $\mathcal{G}_\ell = (\mathcal{D}_\ell, E_\ell)$. An undirected edge exists between $Z_{i_1,j_1}^\ell$ and $Z_{i_2,j_2}^\ell$ if and only if they form a partial silting pair, evaluated via \Cref{thm:pairwise}\textup{(A, B1)}.
\item \textbf{Maximal clique enumeration.} Run the Bron--Kerbosch algorithm on the neighbourhood of $Z_{0,0}^\ell$ in $\mathcal{G}_\ell$ to find all maximal cliques containing $Z_{0,0}^\ell$, and filter the output to retain only cliques of size $n+m$.
\item \textbf{Canonicalisation.} Quotient the output by horizontal lattice translation $(i,j) \mapsto (i+c, j)$, translating each clique so its maximum $i$-coordinate at $j=0$ equals $0$. (Note: horizontal translation is a combinatorial symmetry of the pairwise criterion, not a categorical $\Sigma$-equivalence).
\end{enumerate}

\begin{proposition}[Complexity analysis]\label{prop:complexity}
For fixed $\ell$, Algorithm B.1 operates on a graph with
\[
|\mathcal D_\ell| = b\,(2a+1)
\]
vertices, and its worst-case running time is $O\bigl(3^{\,b(2a+1)/3}\bigr)$.
\end{proposition}

\begin{proof}
The vertex count is immediate from Step 2. The Bron--Kerbosch algorithm finds all maximal cliques in a graph on $N$ vertices in worst-case time $O(3^{N/3})$. Substituting $N = b(2a+1)$ gives the stated bound. The shift in Step 5 quotients out $\Sigma$-equivalence, so the output is a set of representatives of maximal intra-component partial silting configurations.
\end{proof}

\begin{remark}
The algorithm is deliberately restricted to a single component, because the domain restriction in Step 2 is derived from the intra-component staircase bound. Extending it to configurations spread across several components requires replacing $\mathcal D_k$ with a domain cut out by the cross-component inequalities of \Cref{thm:cross-pairwise}; the resulting vertex count, and hence the complexity bound, then depends on the number of components involved.
\end{remark}





\end{document}